\documentclass [a4paper, twoside]{article}
\usepackage{latexsym}
\title{ 
Improving the primal-dual algorithm for the transportation problem in
the plane
}

\author{ Thomas Kaijser
\\ 
\it\small Department of Mathematics
and
Information Coding Group,
\\
\it\small Link\"{o}ping University, S-581 83 Link\"{o}ping, Sweden
\/{\rm ;} \it\small thkai@mai.liu.se }

\date{}
\begin{document}
\maketitle
\newtheorem{lem}{Lemma}
\newtheorem{thm}{Theorem}
\newtheorem{corr}{Corollary}
\newtheorem{path}{Proposition}

\begin{abstract}
The transportation problem in the plane - how to move a set of objects
from one set of points to another set of points in the cheapest way -
is a very old problem going
back several hundreds of years. In recent years the solution of the
problem has found 
applications in the analysis of digital images when searching for
similarities and discrepancies between images. The main drawback,
however, is the long computation time for finding the solution.

In this paper we present some new results by which the time for
solving the transportation problem
in the plane can be reduced substantially. As cost-function we 
choose a distance-function between points in the plane. We consider
both the case when
the distance-function is equal to the ordinary Euclidean distance, as
well as the case when the distance-function is equal to the square
of the Euclidean distance. This latter distance-function has the
advantage that it
is integer-valued if the coordinates of the points in the plane are 
integers.

\vspace{.5cm}
{\bf Mathematics Subject Classification (2000)}: Primary 90C08; 
Secondary 90C46, 90C06, 90C90.
\[\]


\end{abstract}

\section {Introduction.}

The classical transportation problem can be formulated
as follows. Let \(\{S_{n},\; n=1,2,...,N\}\) denote \(N\) sources, 
let  \(\{R_{m}, \;m=1,2,...,M\}\) denote \(M\) sinks, let 
\(a_{n}\) denote the amount of goods that is available at source
\(S_{n}\), let 
\(b_{m}\) denote the amount of goods required at the sink
\(R_{m}\),  and let \(c(n,m)\) denote the cost to send one unit
of goods from the source 
\(S_{n}\) to the sink
\(R_{m}\). We call \(a_{n}\) the {\em storage} at source 
\(S_{n}\) and
\( b_{m}\) the {\em demand} at sink \(R_{m}\).
Assume that \(\sum_{1}^{N} a_{n} = \sum_{1}^{M} b_{m}\) and let 
\(\Lambda\) 
denote the set of \(N \times M\) matrices 
\(\{y(n,m),\;n=1,2,...,N, \;m=1,2,...,M\}\), 
such that
\(\sum_{m=1}^{M}y(n,m)=a_{n}\), 
\(\sum_{n=1}^{N}y(n,m)=b_{m}\) and \(y(n,m)\geq 0, 
\;n=1,2,...,N, \;m=1,2,...,M.\) 
We call a matrix in 
\(\Lambda\) a {\em complete transportation plan}. 
\newline
Problem: Find a complete transportation
plan  
\(X=\{x(n,m)\}\in\Lambda\) such that
\[
   \sum_{n}\sum_{m}x(n,m)c(n,m) 
   = \min \{\sum_{n}\sum_{m}y(n,m)c(n,m), \;\;\{y(n,m)\}\in \Lambda\}.
\]

This problem, the (balanced) transportation problem, is a basic example 
within the theory of optimization theory, and in standard text books
one usually presents two solution methods to the problem, namely 
the simplex method and the primal-dual algorithm. The simplex method
was developed by G. Danzig in the late 1940s. 
An early presentation of the primal-dual algorithm 
for the special
version of the transportation problem called the {\em assignment problem}, 
which occurs if \(N=M\) and 
\(a_{n}=1\), for \( n=1,2,...,N\) and \(b_{m}=1,\) for \(m=1,2,...,M \), 
was given by 
H. Kuhn in \cite{Kuh55}.
In this paper Kuhn writes "One interesting aspect of the algorithm is
the fact that it is 
latent in work of D. K\"{o}nig and 
E. Egerv\'{a}ry
 that predates the birth of linear programming
by more than 15 years (hence the name, the `{\em Hungarian algorithm}')".

In this paper we shall consider the sources and the sinks 
as points in the plane and write
\(\{S_{n}=(i_{n},j_{n}),n=1,2,...,N\}\) for the sources and 
\(\{R_{m}=(x_{m},y_{m}),m=1,2,...,M\}\) for the sinks.
 We shall also assume that the
cost-function is a distance-function 
\(\delta((i,j),(x,y))\) 
between 
points in the plane. We call this special version of
the transportation problem for
the {\em transportation problem in the plane}. (From now on we write
\(\delta(i,j,x,y)\) 
instead of 
\(\delta((i,j),(x,y))\) for a distance-function between points in the  
plane.)

The transportation problem has a long history, 
and goes as far back as to Monge, 1781 \cite{Mon81}.

In a paper from 1942 L. Kantorovich proves that the solution of the 
transportation problem, can be obtained by solving a maximization problem
instead of a minimization problem - the so called dual version of 
the transportation problem. (See \cite{Kan42}.) He also
showed that the solution to the transportation problem can be
used as a distance-measure between probabilities.

There is much literature on the transportation problem. Here we only mention
the paper \cite{Rac85} and the
book \cite{Rac91} by Rachev in both of which there are many further
references.

In the 1980s the solution to the transportation problem in the plane  
was introduced as a distance-measure for digital, grey-valued images
with equal total grey-value. 
See \cite{WPR85} and \cite{Wer86}. We call this distance-measure
the {\em Kantorovich distance for images}.

A drawback with the Kantorovich distance for
images is, that in case we deal with ordinary grey-valued
images, the size of the transportation problem becomes quite large
which implies that the computation time becomes long.

Standard algorithms for finding the solution to the transportation problem -
namely the simplex method and the primal dual algorithm, both have a 
computational complexity of order \(O(N^{3})\) 
in case the data is integer valued, a fact which was pointed out already by
Werman et al. in \cite{WPR85}. 

In 1995, Atkinson and Vaidya \cite{AV95} presented an
algorithm for the transportation problem in the plane, which, for integer
valued problems, has
a computational complexity of order \(O(N^{2}\times log(N) \times log(C))\), 
where \(C\) denotes the maximum of the storages and the demands, in case
the distance-function
is the \(L^{1}-metric\), and
of order \(O(N^{2.5}\times(log(N))^{3}\times log(C))\) in case 
the distance-function
is {\em the Euclidean metric}. However, they did not
apply their algorithm to digital images, nor did they actually present
any experimental results, so it is difficult to make comparisons with
other methods.

In a recent paper \cite{Kai98} from 1998 (see also \cite{Kai96}), we 
described an algorithm for computing the Kantorovich distance for
images, 
in case the distance-function 
between points in the plane is chosen to be the
\(L^{1}-metric\) or {\em the square of the
Euclidean metric}; the algorithm was 
based on the
primal-dual algorithm for the balanced transportation problem,
and practical experiments indicated a computational complexity of
order approximately \(O(N^{2})\) for square images with \(N\) pixels.

The reason we managed to obtain an algorithm, which 
in comparison with the standard primal-dual algorithm for the
balanced transportation problem,
has a lower computational complexity, 
was because we found a method by which we could reduce the search 
for the
so called \emph{admissible arcs}.

However, it was only in the case
when the distance-function
 is the 
 \(L^{1}-metric\) we were able \emph {to prove} that our algorithm
computed the Kantorovich distance exactly.  In case the 
distance-function 
\(\delta(x_{1},y_{1},x_{2},y_{2})\) is defined by
\begin{equation} 
\delta(x_{1},y_{1},x_{2},y_{2})=(x_{1}-x_{2})^{2}+(y_{1}-y_{2})^{2}
\end{equation} 
we were only able to show that our
computation leads to the correct result by checking an optimality
criterion which exists for the primal-dual algorithm. (As we shall
motivate in Section 6 below, there are several reasons why it is of
interest to use a distance-function defined by (1).)
 
The purpose of this paper is threefold. 
One purpose of this paper is to prove that the stopping criterion
introduced in \cite{Kai98}  
does  indeed lead to the correct result 
in case the distance-function is defined
by (1).

A second purpose is to prove some results by
which one can \emph{improve} the algorithm for computing the Kantorovich 
distance for images in case the distance-function again is defined by (1).
These results imply that we
can restrict the search for admissible arcs even further. 

One case we did not consider in the paper 
\cite{Kai98} 
was the
case when the distance-function is exactly the Euclidean distance.
(We mentioned very briefly, \cite{Kai98} Section 20, 
that in case one uses a linear 
combination, of the \(L^{\infty}-metric\) and 
the \(L^{1}-metric\) then the algorithm we described in
\cite{Kai98} did work for the examples we had considered.)

The third purpose of this paper is to prove 
a result by which it is possible to
improve the primal-dual algorithm for 
computing the Kantorovich distance also in
case the distance-function is the Euclidean distance.
We believe, that by using this result, it should be possible to obtain
an algorithm for the transportation problem in the plane with 
Euclidean distance-function, which, for integer valued problems, has 
a computational complexity of approximately order \(O(N^{2})\) to
be compared to \(O(N^{2.5})\) 
obtained approximately in \cite{AV95}.

The plan of the paper is as follows. In Section 2 we introduce 
some basic notations, some terminology and a more precise 
formulation of the transportation problem in the plane. In Section
3 we present the dual formulation of the transportation problem 
and in Section 4 we give a very
brief description of the primal-dual algorithm. 
In Section 5 we recall some results proven in \cite{Kai98}.
In Section 6 
we prove the assertion which we presented in \cite{Kai98} in case the 
distance-function is defined by (1), and in Section 7 we
prove some other results which can be used
to decrease the search for admissible arcs when the distance-function
is defined by (1). 
In Section 8 we prepare for our results, when the 
distance-function 
is the Euclidean distance, by introducing three notions namely 
a \emph{hyperbolic set}, a \emph{level set} and an 
\emph{exclusion set}, and we prove a simple but very important
relation between hyperbolic sets and level sets. In Section 9, 
we 
prove
a simple, elementary, lemma  which relates hyperbolas and cones, and then 
in Section 10, we apply this result when proving results, by which one can
reduce the search for admissible arcs in case the 
distance-function is the ordinary Euclidean distance.
In Section 11, for the sake of completeness,  
we prove that one of the basic assumptions we make in most
of our results is correct, namely that if one uses the primal-dual
algorithm with proper initialization then each pixel will always 
belong to at least 
one admissible arc. 
In section 12 we present a result 
by which one can speed up the computation of the quantities by which
one changes the dual variables when performing the primal-dual algorithm.

In order for the results of this paper to be useful it is necessary that
the sources and the sinks can be organized in such a way that one quickly
can determine all points which are so to speak northeast (northwest, 
southeast, southwest) of an arbitrary point.   
In section 13 we briefly describe an algorithm by which one can 
accomplish this.
Section 14, finally, contains a short summary.

\section{Basic notations and terminology.}
Let \(K=\{(i_{n},j_{n}), \;
n=1,2,...,N \}\)
 be a set of
\(N\) points in \(R^{2}\). By an \emph{image} \(P\) with support \(K\) 
(defined on \(K\))
we simply 
mean a set \(P=\{(i_{n},j_{n},p(i_{n},j_{n})), \;n=1,2,...,N \}\) for
which \(p(i_{n},j_{n})>0, \;n=1,2,...,N \). 
We usually use the notation \(P=\{p(i,j):\;(i,j)\in K\}\) for an image.
We call an element \((i,j)\) of the support \(K\) of an image a
{\em pixel}.
  
Next, let
\(P=\{p(i,j):\;(i,j)\in K_{1}\}\) 
and
\(Q=\{q(x,y):\;(x,y)\in K_{2}\}\)  
be two given images, defined on the two sets 
\(K_{1}=\{(i_{n},j_{n}), \;n=1,2,...,N\}\)
 and 
\(K_{2}=\{(x_{m},y_{m}), \;m=1,2,...,M\}\)
respectively.
\(K_{1}\) and \(K_{2}\)
may be the same, overlap or be disjoint. 
We also assume that 
\[ \sum_{K_{1}}p(i,j) =\sum_{K_{2}}q(x,y). \]
Let 
\(\Gamma(P,Q)\) 
denote the set 
of all non-negative mappings \(h(i,j,x,y)\) from 
\( K_{1} \times K_{2} \rightarrow R^{+}\) such that 
\begin{equation}
\sum_{m=1}^{M} h(i(n),j(n),x(m),y(m))\leq p(i(n),j(n)), \;\;n=1,2,...,N
\end{equation}
and  
\begin{equation}
\sum_{n=1}^{N} h(i(n),j(n),x(m),y(m))\leq q(x(m),y(m)) , \;\;m=1,2,...,M.
\end{equation}
We call any function in 
\(\Gamma(P,Q)\) a {\em transportation plan} from \(P\) to \(Q\). 
A transportation plan
for which we have equality in both (2) and (3) will be called 
a {\em complete transportation plan} and we denote the set of
all complete transportation plans by 
\(\Lambda(P,Q)\).

Let
\(\delta(i,j,x,y)\) be a {\em distance-function} 
measuring the distance from a pixel \((i,j)\) in \(K_{1}\) 
to a pixel \((x,y)\) in 
\(K_{2}\). 
The {\em Kantorovich distance}
\(d_{K,\delta}(P,Q)\) with underlying distance-function
\(\delta(i,j,x,y)\)
is defined by 
\[ d_{K,\delta}(P,Q) = 
\min \{\sum_{i,j,x,y}
 h(i,j,x,y)\times \delta(i,j,x,y):
\;\; 
h(\cdot,\cdot,\cdot,\cdot) \in \Lambda(P,Q)\}.\]
From the definition and some consideration we see
 that computing the Kantorovich distance
for images is equivalent to solving a linear programming problem 
called the {\em balanced transportation
problem}. Most standard text books in optimization theory presents
an algorithm for solving the general balanced transportation problem
with arbitrary cost-function. See e.g. \cite{AMO93}, \cite{Ber91} or
\cite{Mur76}.

\section{The dual formulation.}
  It is well-known that the solution to 
the minimization problem described above is also obtained by solving
the following maximization problem - the dual problem:
\[ maximize \;\;\sum_{(i,j)\in K_{1}}
 \alpha(i,j)\times p(i,j) + 
\sum_{(x,y)\in K_{2}}
 \beta(x,y)\times q(x,y) \]
\(when\)
\begin{equation}\delta(i,j,x,y) - \alpha(i,j)- \beta(x,y)\ge 0, \;\; \;\;
(i,j) \in K_{1}, \; (x,y) \in K_{2}. \end{equation}
The variables 
\( \alpha(i,j)\)
and
 \(\beta(x,y)\) 
are called the {\em dual} variables and in case they
satisfy (4) we have a dual {\em feasible} solution.
A pair \(\{(i,j),(x,y)\}\) where 
\((i,j) \in K_{1}\) and \((x,y) \in K_{2}\) is called an {\em arc}.
In case  
\begin{equation}
\alpha(i,j)+ \beta(x,y)= \delta(i,j,x,y)
\end{equation}
we say that 
\(\{(i,j),(x,y)\}\)
 is an {\em admissible} arc. Otherwise the arc is
called {\em unadmissible}. 

If a pixel \((i,j) \in K_{1}\) is such that there exists 
a pixel \((x,y) \in K_{2}\) such that (5) holds then we say that
\((i,j)\) {\em has } an admissible arc and vice versa.
If  
\(\{(i,j),(x,y)\}\) 
is an {\em admissible} arc we also say that
\((x,y)\) is an {\em admissible pixel} with respect to 
\((i,j)\). And vice versa.
\section{The primal-dual algorithm.}
The primal-dual algorithm runs roughly as follows. Let 
\(\{\alpha(i,j), \beta(x,y): \;(i,j) \in K_{1}, \; (x,y) \in K_{2}\}\)
be a feasible set of dual variables, and let us also assume that initially 
there 
exists at least one admissible arc for every pixel in 
\(K_{1}\) and  \(K_{2}\). (If this is not the case originally, it is easy to 
see that one can increase some or all of the dual variables
so that this hypothesis is
fulfilled). We now look for a transportation plan \(h(i,j,x,y)\) 
from \(P\) to \(Q\) 
which has the largest total mass among all the transportation plans
for which \(h(i,j,x,y) = 0\) in case \(\{(i,j),(x,y)\}\) is
unadmissible. 
In case the transportation plan we find is complete then we are
ready; otherwise we update our set of dual variables. 
To find the "best" 
transportation plan on a given set of admissible arcs we use a labeling 
process, and we use the labeling also in order to find the quantity
used when updating the set of dual 
variables. 
Once we have updated the variables, we determine the
new set of admissible arcs and then we again look 
for the "best" transportation plan on 
this new set of admissible arcs. Etcetera.

A good reference on the primal-dual algorithm for the transportation 
problem is Murty \cite{Mur76}, chapter 12. 
Other useful references are \cite{Ber91} and 
\cite{AMO93}.
\newline 
\[\]
\section{Some preliminary results.} 
In \cite{Kai98} we were able to 
improve the primal-dual algorithm for computing the Kantorovich
distance for images. The reason for the improvement was that we
were
able to reduce the search for new admissible arcs,
by using the structure of the underlying distance-functions.
In this section
we shall present some results and proofs which in essence can be
found in \cite{Kai98}. 

The following proposition holds for any choice of underlying
distance-function.

\begin{path} Let the underlying distance-function \(\delta(i,j,x,y)\) be
arbitrary. 
Let \((i_{1},j_{1})\) and 
\((i_{2},j_{2})\) belong to 
 \(K_{1}\), let \((x_{1},y_{1})\)
and 
\((x_{2},y_{2})\) belong to
 \( K_{2}\) and suppose that 
\(\{(i_{2},j_{2}),(x_{2},y_{2})\}\) 
is an admissible arc.
Then 
\begin{equation}
    \alpha(i_{1},j_{1})- \alpha(i_{2},j_{2})
     \le \delta(i_{1},j_{1},x_{2},y_{2})-
\delta(i_{2},j_{2},x_{2},y_{2})
\end{equation}
and similarly
\begin{equation}
\beta(x_{1},y_{1})- \beta(x_{2},y_{2})
\le \delta(i_{2},j_{2},x_{1},y_{1})-
    \delta(i_{2},j_{2},x_{2},y_{2}).
\end{equation}
\end{path}
{\bf Proof.}
Let us prove (7). 
We have
\[\beta(x_{1},y_{1})- \beta(x_{2},y_{2}) = \]
\[\beta(x_{1},y_{1})- \delta(i_{2},j_{2},x_{2},y_{2}) + 
\alpha(i_{2},j_{2}) 
\leq \]
\[
\delta(i_{2},j_{2},x_{1},y_{1}) - \alpha(i_{2},j_{2}) -
\delta(i_{2},j_{2},x_{2},y_{2}) +
\alpha(i_{2},j_{2}) = \]
\[
\delta(i_{2},j_{2},x_{1},y_{1}) - 
\delta(i_{2},j_{2},x_{2},y_{2}).\]
The proof of (6) can be done in an analogous way. QED.

By applying the triangle inequality the following proposition 
follows immediately from
Proposition 1. We therefore
state it without proof.

\begin{path} Let the underlying distance-function \(\delta(i,j,x,y)\) be
a metric. 
Then, if both 
\((i_{1},j_{1}) \in K_{1}\) 
and  
\((i_{2},j_{2})\in K_{1}\) have admissible arcs then 
\begin{equation}
    |\alpha(i_{1},j_{1})- \alpha(i_{2},j_{2})|
     \le \delta(i_{1},j_{1},i_{2},j_{2}),
\end{equation}
and similarly, if both 
\((x_{1},y_{1}) \in K_{2}\) 
and 
\((x_{2},y_{2})\in K_{2},\) 
have admissible arcs 
then
\begin{equation}
    |\beta(x_{1},y_{1})- \beta(x_{2},y_{2})|
     \le \delta(x_{1},y_{1},x_{2},y_{2}).
\end{equation}
\end{path}

Before we state and prove the next lemma let us introduce some
further terminology. 

Let \((i,j)\) be a pixel in the support \(K_{1}\)
of the image \(P\) and let \(\alpha(i,j)\) be a 
dual variable corresponding to the pixel \((i,j)\). 
Let 
\((x,y)\) 
be a pixel in the support of the image \(Q\).
If the dual variable \(\beta(x,y)\) is such that
\[ \beta(x,y) < \delta(i,j,x,y)-\alpha(i,j)\]
then we say that 
\((x,y)\) is {\em low 
with respect to} \((i,j)\). 
In case there is little risk for misunderstanding we only say that
\((x,y)\) is {\em low}. 

Let us also introduce the following notation and terminology regarding the
positions of two pixels. Thus let
\((x_{1},y_{1})\)
and 
\((x_{2},y_{2})\)
be two pixels. If  
\(x_{1}\le x_{2}\) and \(y_{1}\le y_{2}\) then
we say that  
\((x_{2},y_{2})\)
is {\em northeast (NE)} of
\((x_{1},y_{1})\), and that 
\((x_{1},y_{1})\) is {\em southwest (SW)} of
\((x_{2},y_{2})\), and if 
\(x_{1}\ge x_{2}\) and \(y_{1}\le y_{2}\) then
we say that  
\((x_{2},y_{2})\)
is {\em northwest (NW)} of
\((x_{1},y_{1})\), and that 
\((x_{1},y_{1})\) is {\em southeast (SE)} of
\((x_{2},y_{2})\).

The usefulness of our next result is that it helps to limit the 
number of tests needed for finding all new admissible arcs 
in case we use the
\(l_{1}\)-metric as underlying distance-function. 
The result can be found in \cite{Kai98}, (Lemma 19.1)  but we repeat it here
as background information.

\begin{thm}
Suppose that the distance-function we are using is 
defined by 
\[ \delta(i,j,x,y)=|i-x|+|j-y|.\]
Let (i,j) be a pixel in \(K_{1}\), 
which has an admissible arc.
Now suppose that 
\((x_{1},y_{1}) \in K_{2}\),
that  \((x_{1},y_{1})\) has an admissible arc ,
that \((x_{1},y_{1})\) is NE of \((i,j)\) and 
that \((x_{1},y_{1})\) 
is low with respect to 
\((i,j)\). Then, if \((x,y)\) is NE of 
 \((x_{1},y_{1})\), then  
\((x,y)\) is also low with respect to
\((i,j)\). (See Figure 1 below.) 
\end{thm}

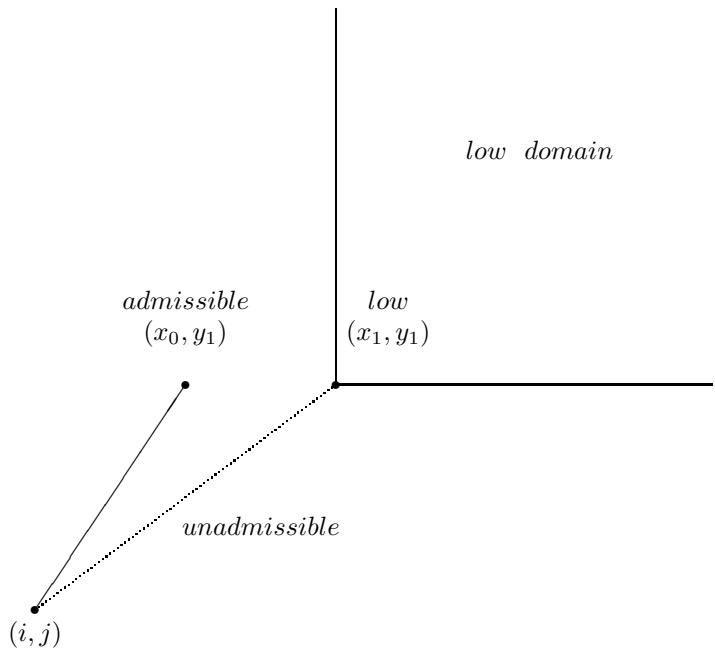
{\bf Figure 1}

\setlength{\unitlength}{1mm}
\begin{picture}(100,100)(-10,-10)
\linethickness{1pt}
\thinlines
\put(0,0){\circle*{1}}
\put(0,-5){\makebox(0,0)[b]{\((i,j)\)}}
\put(20,30){\circle*{1}}
\put(20,35){\makebox(0,0)[b]{\((x_{0},y_{1})\)}}
\put(20,40){\makebox(0,0)[b]{\(admissible\)}}
\put(30,10){\makebox(0,0)[b]{\(unadmissible\)}}
\put(40,30){\circle*{1}}
\put(47,35){\makebox(0,0)[b]{\((x_{1},y_{1})\)}}
\put(47,40){\makebox(0,0)[b]{\(low\)}}
\put(60,60){\makebox(0,0)[b]{\(low \)}}
\put(71,60){\makebox(0,0)[b]{\(domain \)}}
\put(0,0){\line(2,3){20}}
\qbezier[80](0,0)(20,15)(40,30)
\put(40,30){\line(1,0){50}}
\put(40,30){\line(0,1){50}}

\end{picture}

{\bf Proof.} 
We prove 
the theorem
by contradiction. Thus suppose that there exists
a pixel \((x,y)\in K_{2}\) located \(NE\) of 
\( (x_{1},y_{1})\) 
and
such that at that pixel 
the dual variable \(\beta(x,y)\) is such that
\begin{equation}
-\delta(i,j,x,y) + \alpha(i,j)+ \beta(x,y) = 0. 
\end{equation}
But since \( (x_{1},y_{1})\) is low with respect to \((i,j)\), it
follows that 
\[\alpha(i,j) < \delta(i,j,x_{1},y_{1}) -\beta(x_{1},y_{1}) \]
which together 
with 
(10) implies that 
\[\delta(i,j,x,y) - \beta(x,y) 
< \delta(i,j,x_{1},y_{1}) -\beta(x_{1},y_{1}) \]
and hence
\[ \beta(x,y) - \beta(x_{1},y_{1}) >
\delta(i,j,x,y) - \delta(i,j,x_{1},y_{1}).\] But since \((x,y)\)
is \(NE\) of 
\((x_{1},y_{1})\) which is \(NE\) of \((i,j)\) it follows that
\[\delta(i,j,x,y) - \delta(i,j,x_{1},y_{1})
=
x-i+y-j-
(x_{1}-i+y_{1}-j)=\]
\[
x-x_{1}+y-y_{1}=\delta(x,y,x_{1},y_{1}).\]
Hence
\[ \beta(x,y) - \beta(x_{1},y_{1}) > \delta(x,y,x_{1},y_{1}).\]
But since both \((x,y)\) and 
\((x_{1},y_{1})\) have admissible arcs this inequality can not
be fulfilled because of Proposition 2.
QED.

\section {Finding admissible arcs when the distance-function 
is the square of the Euclidean distance.}

In
\cite{Kai98} we also considered the case when 
the underlying distance-function is defined as 
{\em the square of the Euclidean distance}. 
In our computer experiments we first tried to apply the assertion of 
Theorem 1 in order
to reduce the search for new admissible arcs, but it turned out, that
if we did so, then, for some examples, we did in fact exclude too many arcs.
Therefore we created another
claim with slightly stronger assumptions, a claim which we formulated
as an Assertion, (\cite{Kai98} Assertion 20.1).
It turned out that by using this assertion
in order to reduce the search for admissible arcs,
we did not exclude any admissible
arcs from our search. 
However we were at that time not able to prove the claim, but now we have
a proof and can therefore formulate our claim as a theorem.
We shall first introduce yet 
some further
terminology.

Consider two pixels
\((x_{1},y_{1})\) 
and
\((x_{2},y_{2})\) in \(K_{2}\). 
Suppose both 
\((x_{1},y_{1})\) 
and
\((x_{2},y_{2})\) 
are low with respect to the pixel 
\((i,j)\in K_{1}\). If also 
\[\delta(i,j,x_{2},y_{2}) - \alpha(i,j)- \beta(x_{2},y_{2}) \geq \]
\begin{equation}
\delta(i,j,x_{1},y_{1}) - \alpha(i,j)- \beta(x_{1},y_{1})  
\end {equation}
then we say that 
\((x_{2},y_{2})\) 
is {\em  lower} than 
\((x_{1},y_{1})\). If we have strict inequality in (11) we  
say that 
\((x_{2},y_{2})\) 
is {\em  strictly lower} than 
\((x_{1},y_{1})\).

\begin{thm}
Suppose that the underlying distance-function is defined by
\begin{equation}
\delta(i,j,x,y)=(i-x)^{2}+(j-y)^{2}.
\end{equation}
Let 
\(\{\alpha(i,j), \beta(x,y)\}\) 
be a feasible set of dual variables such that 
each pixel \((i,j)\) in \(K_{1}\), and 
each pixel \((x,y)\) in \(K_{2}\) have an admissible arc.
Let \((i_{0},j_{0})\) be a pixel in \(K_{1}\) and let
\((x_{1},y_{1})\),
\((x_{2},y_{1})\) and
\((x_{3},y_{1})\) 
be three pixels in  \(K_{2}\) on the same line \(y=y_{1}\), such
that \(i_{0}\leq x_{1} < x_{2} < x_{3}\). Furthermore assume that
both 
\((x_{1},y_{1})\) and 
\((x_{2},y_{1})\) are low with respect to \((i_{0},j_{0})\) and let
\((x_{2},y_{1})\) be lower than 
\((x_{1},y_{1})\).

Then, 
\((x_{3},y_{1})\) is also low with respect to \((i_{0},j_{0})\). (See
Figure 2 below, for a graphical illustration.)
\end{thm}
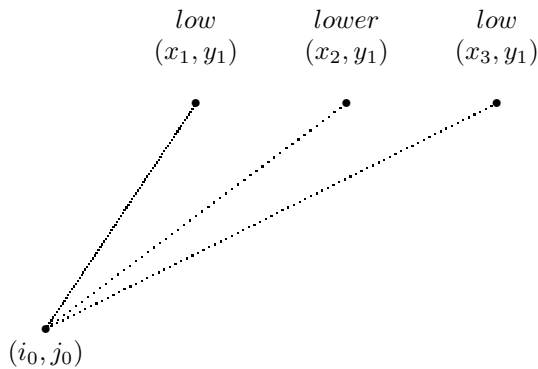
{\bf Figure 2.}
\setlength{\unitlength}{1mm}
\begin{picture}(65,65)(-10,-10)
\linethickness{1pt}
\thinlines
\put(0,0){\circle*{1}}
\put(0,-5){\makebox(0,0)[b]{\((i_{0},j_{0})\)}}
\put(20,30){\circle*{1}}
\put(20,35){\makebox(0,0)[b]{\((x_{1},y_{1})\)}}
\put(20,40){\makebox(0,0)[b]{\(low\)}}
\put(40,30){\circle*{1}}
\put(40,35){\makebox(0,0)[b]{\((x_{2},y_{1})\)}}
\put(40,40){\makebox(0,0)[b]{\(lower\)}}
\put(60,30){\circle*{1}}
\put(60,35){\makebox(0,0)[b]{\((x_{3},y_{1})\)}}
\put(60,40){\makebox(0,0)[b]{\(low\)}}
\qbezier[80](0,0)(10,15)(20,30)
\qbezier[50](0,0)(20,15)(40,30)
\qbezier[80](0,0)(30,15)(60,30)
\end{picture}
\[\]
{\bf Remark 1.} This theorem is slightly sharper than what we
formulate in Assertion 20.1 of \cite{Kai98} ( see also \cite{Kai96},
Condition 18.1), since we only require
that \((x_{2},y_{1})\) shall be lower
than \((x_{1},y_{1})\) and not necessarily strictly lower. 
\newline
{\bf Remark 2.} This theorem requires that our images have supports on
a regular grid structure where it is meaningful to speak of pixels 
being located on the same line. The theorem is not particularly useful 
for the 
case when the elements of the supports of the two images under
consideration are placed at random. 
\newline
{\bf Remark 3.} 
It ought to be clear that by using this theorem one can decrease
 the number of
searches necessary to find all admissible arcs. On each line 
\(y=y_{1}\)
we do not have to check pixels "further away", once we have found
two pixels \((x_{1},y_{1})\),
\((x_{2},y_{1})\) such that both are low and 
\((x_{2},y_{1})\) is lower than \((x_{1},y_{1})\). 
\newline
{\bf Remark 4.} In Section 7 we
shall show how one can obtain an efficient "stopping
criterion" also in the "\(y\)-direction".
\newline
{\bf Remark 5.} There are several reasons for defining the underlying 
distance-function by (5). It is rotationally invariant, 
it takes integers into integers, 
it also gives rise to a metric 
(if one takes the square root afterwards, see e.g. 
\cite{Rac85}), and if the two images are pure translations of each other
then in the general case there is a 
{\em unique} transportation plan which gives
rise to the optimum value.
Moreover, by comparison with the \(L^{1}-metric\) 
it seems to give rise to a distance which is less course.
\[\] 
Before we begin our proof of Theorem 2 we shall prove the following
auxiliary result.
\begin{path} a)
Let \((i_{1},j_{1})\) and \((i_{2},j_{2})\) be two pixels 
in \(K_{1}\), let
\((x_{1},y_{1})\) and 
\((x_{2},y_{1})\) be two pixels in
\(K_{2}\) on the same line and such that 
\(x_{1} < x_{2}\).
Let 
\(\delta(\cdot,\cdot,\cdot,\cdot)\) 
be defined by (12)
and suppose 
that \(\{(i_{1},j_{1}),(x_{1},y_{1})\}\) and  
\(\{(i_{2},j_{2}),(x_{2},y_{1})\}\) are admissible arcs.
Then
\(i_{1} \leq i_{2}\).
\newline
b) If furthermore  
\(\{(i_{2},j_{2}),(x_{1},y_{1})\}\) is not an admissible arc
then
\(i_{1} < i_{2}\).
\newline
(For a graphical illustration, see Figure 3 below.)

\end{path}
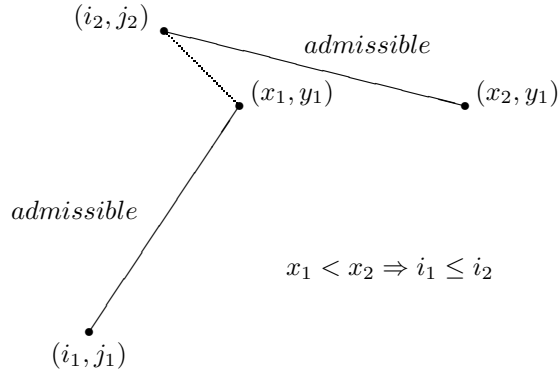
{\bf Figure 3}

\begin{picture}(70,70)(-20,-20)
\linethickness{1pt}
\thinlines
\put(0,0){\circle*{1}}
\put(0,-5){\makebox(0,0)[b]{\((i_{1},j_{1})\)}}
\put(20,30){\circle*{1}}
\put(27,30){\makebox(0,0)[b]{\((x_{1},y_{1})\)}}
\put(50,30){\circle*{1}}
\put(57,30){\makebox(0,0)[b]{\((x_{2},y_{1})\)}}
\put(0,0){\line(2,3){20}}
\put(10,40){\circle*{1}}
\put(3,40){\makebox(0,0)[b]{\((i_{2},j_{2})\)}}
\put(37,37){\makebox(0,0)[b]{\(admissible\)}}
\put(-2,15){\makebox(0,0)[b]{\(admissible\)}}
\put(10,40){\line(4,-1){40}}

\qbezier[30](10,40)(15,35)(20,30)
\put(40,7){\makebox(0,0)[b]{\(x_{1}<x_{2}\Rightarrow i_{1}\leq i_{2}\)}}
\end{picture}

{\bf Proof of Proposition 3.} We first prove part a).
The following relations must hold: 
\[
\alpha(i_{1},j_{1})+\beta(x_{1},y_{1})=
\delta(i_{1},j_{1},x_{1},y_{1})
\]
\[\alpha(i_{2},j_{2})+\beta(x_{2},y_{1})=
\delta(i_{2},j_{2},x_{2},y_{1})
\]
\[
\alpha(i_{1},j_{1})+\beta(x_{2},y_{1}) \leq
\delta(i_{1},j_{1},x_{2},y_{1})
\]
\begin{equation} 
\alpha(i_{2},j_{2})+\beta(x_{1},y_{1}) \leq
\delta(i_{2},j_{2},x_{1},y_{1}).
\end{equation}  
By first adding the last two inequalities and then subtracting the first
two equalities, and finally shifting terms we obtain the inequality
\begin{equation} 
\delta(i_{1},j_{1},x_{2},y_{1})-
\delta(i_{1},j_{1},x_{1},y_{1}) \geq
\delta(i_{2},j_{2},x_{2},y_{1})-
\delta(i_{2},j_{2},x_{1},y_{1}).
\end{equation}  
By using the definition of \(\delta(i,j,x,y)\) (see (12)), we find that the
left hand side of (14) becomes equal to 
\[(x_{2}-i_{1})^{2} - (x_{1}-i_{1})^{2} = x_{2}^{2}
-x_{1}^{2} - 2(x_{2}-x_{1}) i_{1}\]
and the right hand side becomes equal to
\[(x_{2}-i_{2})^{2} - (x_{1}-i_{2})^{2} = x_{2}^{2}
-x_{1}^{2} - 2(x_{2}-x_{1}) i_{2}.\]
Hence in order for (14) to be true it is necessary that
\(i_{2} \geq i_{1}\). Thereby part a) is proved.

If furthermore we know that 
\(\{(i_{2},j_{2}),(x_{1},y_{1})\}\) is not an admissible arc, 
then we have strict inequality in (13), which implies that 
we also have strict inequality in (14), and therefore, by using
the same kind of arguments as above, we find that in order
for (14) to be true with strict inequality, it is necessary
that 
\(i_{2} > i_{1}\). Thereby part b) is also proved.
Q.E.D.
\begin{corr} Suppose that the underlying distance-function is defined by (12).
Let 
\(\{\alpha(i,j), \beta(x,y)\}\) 
be a feasible set of dual variables such that 
each pixel \((i,j)\) in \(K_{1}\), and 
each pixel \((x,y)\) in \(K_{2}\) have an admissible arc.
Let \((i_{0},j_{0})\) belong to \(K_{1}\), let 
\((x_{1},y_{1})\),
\((x_{2},y_{1})\),
\((x_{3},y_{1})\)
be three pixels in  \(K_{2}\) on the same line \(y=y_{1}\), such
that \(i_{0}\leq x_{1} < x_{2} < x_{3}\). Furthermore assume that
\((x_{1},y_{1})\) is admissible with respect to
\((i_{0},j_{0})\)
and that 
\((x_{2},y_{1})\) is low with respect to 
\((i_{0},j_{0})\).
Then 
\((x_{3},y_{1})\) is also low with respect to
\((i_{0},j_{0})\).
\end{corr} 
{\bf Proof of Corollary 1.} Assume that
\(\{(i_{0},j_{0}),(x_{3},y_{1})\}\) is an admissible arc.
Let \((i_{2},j_{2})\) be an admissible pixel with respect to
\((x_{2},y_{1})\).
Since 
\((x_{2},y_{1})\) is low with respect to
\((i_{0},j_{0})\), it follows by applying part b) of Proposition 3 to the 
pixels 
\((i_{0},j_{0}),(i_{2},j_{2})\),
\((x_{2},y_{1})\) and  \((x_{3},y_{1})\),
that
\(i_{2}<i_{0}\). On the other hand, by applying part a) of Proposition 3
to the pixels 
\((i_{0},j_{0}),(i_{2},j_{2})\),
\((x_{2},y_{1})\) and  \((x_{1},y_{1})\), we conclude that
\(i_{0}\leq i_{2}\) and thereby we have reached a contradiction. Q.E.D.
\newline
{\bf Remark 1.} The corollary implies  
that in case the pixel 
\((x_{0},y_{1})\in K_{2}\) 
is the first
pixel in \(K_{2}\) 
along the line \(y=y_{1}\) for which \(x_{0} \ge i_{0}\) and 
\((x_{0},y_{1})\)
is an \emph{admissible} pixel with respect to \((i_{0},j_{0})\), 
then we can stop the search for admissible pixels on this
line as soon as we find 
a pixel which is low. If instead 
it turns out that the first pixel
\((x_{0},y_{1})\in K_{2}\), with \(x_{0} \ge i_{0}\),
along the line \(y=y_{1}\), is \emph{not} an admissible pixel with respect 
\((i_{0},j_{0})\), then we should use Theorem 2 until
we find an admissible pixel, after which we can use the corollary. Clearly
the use of the corollary  only makes a minor improvement of the algorithm.
\newline
{\bf Remark 2.} Just as is the case with Theorem 2 the corollary is 
only useful for
images defined on a grid structure.
\[\]
{\bf Proof of Theorem 2.} We are now ready to prove Theorem 2.
Assume that 
\((x_{3},y_{1})\) is admissible with respect to  
\((i_{0},j_{0})\). Let 
\((i_{1},j_{1})\) be an admissible arc with respect to
\((x_{1},y_{1})\) and let 
\((i_{2},j_{2})\) be an admissible arc with respect to
\((x_{2},y_{1})\). 
From a) of Proposition 3 it follows that
\(i_{1}\leq i_{2}\), 
from b) of Proposition 3 we obtain
\(i_{2}< i_{0}\), 
and from the hypothesis of the theorem
we also know that 
\(
i_{0} \leq x_{1} < x_{2} < x_{3}.
\)
Thus 
\begin{equation}
i_{1} \leq i_{2} < i_{0} \leq x_{1} < x_{2} < x_{3}.
\end {equation}  
From the hypothesis of the theorem we also have that:
\begin{equation}
\alpha(i_{0},j_{0}) + \beta(x_{1},y_{1}) =
\delta(i_{0},j_{0},x_{1},y_{1}) - r_{1} 
\end{equation}
\begin{equation}
\alpha(i_{0},j_{0}) + \beta(x_{2},y_{1}) =
\delta(i_{0},j_{0},x_{2},y_{1}) - r_{2} 
\end{equation}
where 
\begin{equation}
r_{2} \geq r_{1} > 0. 
\end{equation}
By subtracting (17) from (16) and using the definition of
\(\delta(i,j,x,y)\) we obtain
\begin{equation}
\beta(x_{1},y_{1}) - \beta(x_{2},y_{1}) =
r_{2} - r_{1} + 2(x_{2}-x_{1}) i_{0} + x_{1}^{2}- x_{2}^{2}\;, 
\end{equation}
and by subtracting the equality
\(
\alpha(i_{2},j_{2}) + \beta(x_{2},y_{1}) =
\delta(i_{2},j_{2},x_{2},y_{1})  
\)
from the inequality
\(
\alpha(i_{2},j_{2}) + \beta(x_{1},y_{1}) \leq
\delta(i_{2},j_{2},x_{1},y_{1})  
\)
we obtain
\begin{equation}
\beta(x_{1},y_{1}) - \beta(x_{2},y_{1}) \leq
2(x_{2}-x_{1}) i_{2} + x_{1}^{2}- x_{2}^{2}\;. 
\end{equation}
Finally, by using (19) to replace the left hand side of (20), we
obtain
\[
r_{2} - r_{1} + 2(x_{2}-x_{1}) i_{0}
\leq  2(x_{2}-x_{1}) i_{2}
\]
which implies that
\[
r_{1} - r_{2} \geq 2(x_{2}-x_{1})(i_{0}-i_{2}) > 0
\]
where the last inequality follows from (15) and hence 
\[
r_{1}> r_{2} 
\]
which violates (18),
and thereby we have reached a contradiction. Q.E.D.

\section{ A stopping criterion along the vertical axis.}  
Lemma 2 and Corollary 1 imply, roughly speaking, that the search time for 
finding new
admissible arcs is decreased from the order 
\(O(N^{2})\) 
to the order
\(O(N^{1.5})\), in case our images are square digital images, 
where \(N\) as before denotes the total number of pixels in the two images. 
(Lemma 2  and Corollary 1 give rise to stopping criteria 
along the horizontal lines.) 
We shall now show, how we can reduce the search time to roughly 
\(O(N)\), by giving a stopping criterion along the vertical lines also.

\begin{thm}
Suppose that the underlying distance-function is defined by
\[
\delta(i,j,x,y)=(i-x)^{2}+(j-y)^{2}.
\]  
Let 
\(\{\alpha(i,j), \beta(x,y)\}\) 
be a feasible set of dual variables such that 
each pixel \((i,j)\) in \(K_{1}\), and 
each pixel \((x,y)\) in \(K_{2}\) have an admissible arc.
Let 
\((i_{0},j_{0})\) 
belong to \(K_{1}\). Let
\((x_{1},y_{1})\) be a pixel in \(K_{2}\), NE of \((i_{0},j_{0})\),
which is not an admissible pixel with respect to 
\((i_{0},j_{0})\).
Let instead 
\((i_{1},j_{1})\)
be an admissible pixel with respect to
\((x_{1},y_{1})\), and assume that also 
\((i_{1},j_{1})\)
is NE of \((i_{0},j_{0})\).
Then if 
\((x_{2},y_{2})\)
is a pixel in \(K_{2}\),
NE of 
\((x_{1},y_{1})\),
then 
\((x_{2},y_{2})\) is lower than 
\((x_{1},y_{1})\).
\newline
(For a graphical illustration, see Figure 4 below.)

\end{thm}

\begin{picture}(100,100)(-10,-10)
\linethickness{1pt}
\thinlines
\put(20,30){\line(1,0){50}}
\put(20,30){\line(0,1){50}}
\put(0,0){\circle*{1}}
\put(0,-5){\makebox(0,0)[b]{\((i_{0},j_{0})\)}}
\put(20,30){\circle*{1}}
\qbezier[100](0,0)(20,30)(20,30)
\put(25,25){\makebox(0,0)[b]{\((x_{1},y_{1})\)}}
\put(25,22){\makebox(0,0)[b]{\(low\)}}
\put(19,8){\makebox(0,0)[b]{\(unadmissible\)}}
\put(36,8){\makebox(0,0)[b]{\(arc\)}}
\put(-2,33){\makebox(0,0)[b]{\(admissible\)}}
\put(12,33){\makebox(0,0)[b]{\(arc\)}}

\put(45,50){\makebox(0,0)[b]{\(low \)}}
\put(56,50){\makebox(0,0)[b]{\(domain \)}}
\put(10,40){\circle*{1}}
\put(3,40){\makebox(0,0)[b]{\((i_{1},j_{1})\)}}
\put(10,40){\line(1,-1){10}}

\end{picture}

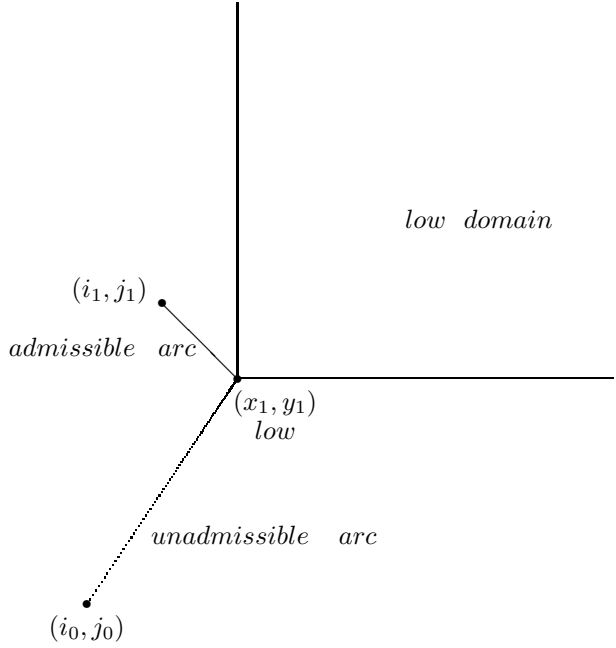
{\bf Figure 4}

\[\]
{\bf Proof.} Let us assume that 
\((x_{2},y_{2})\) 
is \(NE\) of 
\((x_{1},y_{1})\) 
and
that it is 
not lower than 
\((x_{1},y_{1})\).
Then the following relations hold: 
\[
\alpha(i_{0},j_{0})+\beta(x_{1},y_{1})+r_{1}=
\delta(i_{0},j_{0},x_{1},y_{1})
\]
\[
\alpha(i_{1},j_{1})+\beta(x_{2},y_{2}) \leq
\delta(i_{1},j_{1},x_{2},y_{2})
\]
\[\alpha(i_{0},j_{0})+\beta(x_{2},y_{2}) +r_{2}= 
\delta(i_{0},j_{0},x_{2},y_{2})
\]
\[
\alpha(i_{1},j_{1})+\beta(x_{1},y_{1}) =
\delta(i_{1},j_{1},x_{1},y_{1})
\]
and
\[ 0\leq r_{2}<r_{1}.\]
By first adding the first equality and the inequality, and then subtracting 
the next
two equalities, we obtain 
\[
r_{1}-r_{2}\leq \delta(i_{0},j_{0},x_{1},y_{1})
+\delta(i_{1},j_{1},x_{2},y_{2})
-\delta(i_{0},j_{0},x_{2},y_{2})
-\delta(i_{1},j_{1},x_{1},y_{1}),
\]
and shifting terms we find that
\begin{equation}\delta(i_{0},j_{0},x_{1},y_{1})-
\delta(i_{0},j_{0},x_{2},y_{2})
\geq
\delta(i_{1},j_{1},x_{1},y_{1})
-\delta(i_{1},j_{1},x_{2},y_{2})+r_{1}-r_{2}.
\end{equation}
By using the definition of the distance-function 
\(\delta(i,j,x,y)\) 
we find that the left hand side of (21) becomes
equal to
\[
(x_{1}-i_{0})^{2} - (x_{2}-i_{0})^{2} +
(y_{1}-j_{0})^{2} - (y_{2}-j_{0})^{2} 
=\] 
\[
x_{1}^{2}
-x_{2}^{2} - 2(x_{1}-x_{2})i_{0}
+
y_{1}^{2}
-y_{2}^{2} - 2(y_{1}-y_{2}) j_{0},\]
and that the right hand side of (21) becomes equal to 
\[
(x_{1}-i_{1})^{2} - (x_{2}-i_{1})^{2} +
(y_{1}-j_{1})^{2} - (y_{2}-j_{1})^{2} + r_{1}-r_{2}
=
\] 
\[
x_{1}^{2}
-x_{2}^{2} - 2(x_{1}-x_{2}) i_{1}
+
y_{1}^{2}
-y_{2}^{2} - 2(y_{1}-y_{2}) j_{1} + r_{1}-r_{2}.\]
By subtracting the right hand side of (21) from the left hand side,
we thus obtain the inequality
\[
2(x_{1}-x_{2}) i_{1}+2(y_{1}-y_{2}) j_{1}
- 2(x_{1}-x_{2}) i_{0} - 2(y_{1}-y_{2}) j_{0} + r_{2}-r_{1} \geq 0,\]
and by shifting terms we obtain the inequality
\[
2(x_{1}-x_{2}) (i_{1}-i_{0})+2(y_{1}-y_{2}) (j_{1}-j_{0}) \geq r_{1}-r_{2}.
\]
But since \(x_{1}\leq x_{2},\; i_{1}\geq i_{0}, \; 
y_{1} \leq y_{2}\) and \( j_{1} \ge j_{0}\), 
and since \(r_{1}>r_{2}\) by assumption, 
this last inequality can not
be true and we have reached a contradiction. QED.
\newline
{\bf Remark.} This theorem gives rise to a very efficient stopping
criterion since, very often, we can choose the pixel 
\((x_{1},y_{1}) \) equal to \((i_{1},j_{1}) \) equal to 
\((i_{0},y_{1}) \) 
which implies that we can stop the search for
new admissible pixels on or above the line \(y=y_{1}\).
\newline
\[\]
\section{Hyperbolic sets, level sets and
exclusion sets.}
Before we start to discuss how we can shorten the search time
for new admissible arcs in case the underlying distance-function is
the \emph{Euclidean distance}, we shall introduce some further notions.
As usual, let
\(P=\{p(i,j):\;(i,j)\in K_{1}\}\) 
and
\(Q=\{q(x,y):\;(x,y)\in K_{2}\}\)  
be two images, let 
\(\delta(i,j,x,y)\) be an \emph{arbitrary} distance-function, and 
let
\( \{\alpha(i,j):\;(i,j)\in K_{1}\}\) 
and
 \(\{\beta(x,y):\;(x,y)\in K_{2}\}\)  
be two sets of \emph{feasible} dual variables  
associated to the images 
\(P\)
and
\(Q\)
respectively.

For each \((i,j)\in K_{1}\)
we define the \emph{exclusion set}
\(E[(i,j), K_{2}]\) simply as all pixels in 
\(K_{2}\) 
which are not admissible with respect to \((i,j)\) 
and similarly   
for each
\((x,y)\in K_{2}\)  we define the \emph{exclusion set}
\(E[(x,y),K_{1}]\) simply as all pixels in 
\(K_{1}\)
which are not admissible with respect to \((x,y)\).

Next let us also introduce a notion which we call a {\em level set}
as follows. Let \((i_{0},j_{0})\in K_{1}\) and let 
\(r\) be a number such that
\(r\geq 0\). 
The set 
of all pixels \((x,y)\) in 
\(K_{2}\) such that 
\[\delta(i_{0},j_{0},x,y)-\alpha(i_{0},j_{0})- \beta(x,y) > r\]
will be called a {\em 
level set} and denoted by 
\(L[(i_{0},j_{0}),K_{2},r]\).  
 Note that 
\[L[(i,j),K_{2},0] = E[(i,j), K_{2}].\]

We shall now proceed by introducing a notion which we shall
call a \emph{hyperbolic set}.
Let 
\((x_{1},y_{1})\) and \((x_{2},y_{2}) \) 
be two points in \(R^{2}\) 
and let \(r\) be a real number.
We define the \emph{hyperbolic} set  
\(H[(x_{1},y_{1}), (x_{2},y_{2}), r]\) 
by  
\[H[
(x_{1},y_{1}), 
(x_{2},y_{2}), r]=
\{(x,y)\in R^{2}:
\delta(x,y,x_{2},y_{2}) - \delta(x,y,x_{1},y_{1}) < r
\}. \]
\newline
For a graphical illustration, see Figure 5 below.

\[\]
\setlength{\unitlength}{1mm}
\begin{picture}(40,70)(-30,0)
\linethickness{1pt}
\thinlines
\put(0,0){\line(1,4){8}}
\put(0,0){\circle*{1}}
\put(0,-5){\makebox(0,0)[b]{\((x_{1},y_{1})\)}}

\put(10,60){\makebox(0,0)[b]{\(H[(x_{1},y_{2}),(x_{2},j_{2}),r]\)}}

\put(10,55){\makebox(0,0)[b]{\(r < 0 \)}}
\put(8,32){\circle*{1}}
\put(0,32){\makebox(0,0)[b]{\((x_{2},y_{2})\)}}

\linethickness{1pt}
\qbezier[800](-30,80)(-10,-40)(67,67)
\end{picture}
\[\]
{\bf Figure 5.}

\[\]
{\bf Remark.} The reason we call the set an hyperbolic set is
because if we replace the inequality sign in the definition
by an equality sign, thereby changing the defining inequality to an equality,
this equality would define an ordinary hyperbola 
in case the underlying distance-function 
is the ordinary Euclidean distance.
\newline

We shall next present a theorem  which connects hyperbolic sets and
level sets. 

\begin{thm}
Let \(P=\{p(i,j):\;(i,j)\in K_{1}\}\) 
and
\(Q=\{q(x,y):\;(x,y)\in K_{2}\}\)  
be two images, let 
\(\delta(i,j,x,y)\) be an arbitrary distance-function, and 
let
\( \{\alpha(i,j):\;(i,j)\in K_{1}\}\) 
and
 \(\{\beta(x,y):\;(x,y)\in K_{2}\}\)  
be two sets of feasible dual variables  
associated to the images 
\(P\)
and
\(Q\)
respectively, such that 
each pixel \((i,j)\) in \(K_{1}\), and 
each pixel \((x,y)\) in \(K_{2}\) have an admissible arc.

Let 
\((i_{0},j_{0})\) 
be a pixel in \( K_{1}\), let     
\((x_{1},y_{1}) \) 
be a pixel in
\( K_{2}\)
which is  not an admissible pixel with respect to 
\((i_{0},j_{0})\) ,
let the number \(r\) be defined by
\[
r= 
\delta(i_{0},j_{0},x_{1},y_{1}) -
\alpha(
i_{0},j_{0})- 
\beta(x_{1},y_{1}), \]
let 
\((i_{1},j_{1})\) 
be an admissible pixel with respect to the
pixel
\((x_{1},y_{1}), \) 
let \(\Delta\) be defined by
\(\Delta = \delta(i_{0},j_{0},x_{1},y_{1}) -
\delta(i_{1},j_{1},x_{1},y_{1})\),
and let \(s\) be a real number satisfying \(s \geq 0 \).
Then, 
\newline
(a): 
\[H[(i_{0},j_{0}),(i_{1},j_{1}),
r-\Delta- s]
\cap K_{2} \subset
L[(i_{0},j_{0}),K_{2},s],\] 
(b): 
\[H[(i_{0},j_{0}), (i_{1},j_{1}),r-\Delta]\cap K_{2} \subset
E[(i_{0},j_{0}),K_{2}],\] 
and (c):
\newline
if also \(\delta(\cdot,\cdot,\cdot,\cdot)\) is a metric, then
\[H[(i_{0},j_{0}),(x_{1},y_{1}),
r-\delta(i_{0},j_{0},x_{1},y_{1})-s]\subset
H[(i_{0},j_{0}),(i_{1},j_{1}),
r-\Delta-s].\]

\end{thm}
{\bf Proof.} Let us first prove (c). 
We prove (c) by contradiction.
Thus assume that \((x,y)\) does belong to 
\(H[(i_{0},j_{0}),(x_{1},y_{1}),
r-\delta(i_{0},j_{0},x_{1},y_{1})-s]\)
but does not belong to
\(H[(i_{0},j_{0}),(i_{1},j_{1}),
r-\Delta-s]\).
Then
\[
\delta(x,y,x_{1},y_{1})-
\delta(x,y,i_{0},j_{0}) < 
r-\delta(i_{0},j_{0},x_{1},y_{1})-s\]
and also 
\[
\delta(x,y,i_{1},j_{1})-
\delta(x,y,i_{0},j_{0})\geq 
r-\delta(i_{0},j_{0},x_{1},y_{1})+
\delta(i_{1},j_{1},x_{1},y_{1})-s.\]
Hence
\[r-s-\delta(i_{0},j_{0},x_{1},y_{1})+
\delta(i_{1},j_{1},x_{1},y_{1})\leq
\delta(x,y,i_{1},j_{1})-
\delta(x,y,i_{0},j_{0})
=
\]
\[
\delta(x,y,i_{1},j_{1})-
\delta(x,y,x_{1},y_{1}) +
\delta(x,y,x_{1},y_{1}) -
\delta(x,y,i_{0},j_{0})
< \]
\[
r-\delta(i_{0},j_{0},x_{1},y_{1})-s+
\delta(x,y,i_{1},j_{1})-
\delta(x,y,x_{1},y_{1})\]
from which follows after cancellations that
\[
\delta(i_{1},j_{1},x_{1},y_{1})<
\delta(x,y,i_{1},j_{1})-
\delta(x,y,x_{1},y_{1}).
\]
Since we have assumed that 
\(\delta(\cdot,\cdot,\cdot,\cdot)\) is a metric, this last inequality
 can not be true
because of the triangle inequality.

In order to  prove (a)  
let \((x,y)\) denote a pixel in \(K_{2}\) and define
\[r_{2}=
\delta(i_{0},j_{0},x,y)-
\alpha(i_{0},j_{0})-\beta(x,y).
\]
Now suppose that the
pixel \((x,y)\) does not belong to the 
level set
\(L[i_{0},j_{0},K_{2},s]\). This implies that
the number \(r_{2}\) introduced 
above, must satisfy \(0\leq r_{2} \leq s\).
We therefore have the following relations:     

\[\alpha(i_{0},j_{0})+\beta(x,y)+r_{2}=
\delta(i_{0},j_{0},x,y),
\]
\[
\alpha(i_{1},j_{1})+\beta(x_{1},y_{1})=
\delta(i_{1},j_{1},x_{1},y_{1}),
\]
\[
\alpha(i_{0},j_{0})+\beta(x_{1},y_{1}) + r=
\delta(i_{0},j_{0},x_{1},y_{1}),
\]
and 
\[
\alpha(i_{1},j_{1})+\beta(x,y) \leq
\delta(i_{1},j_{1},x,y).
\]
By first adding the first two equations, then subtracting the next equation, 
and finally also subtracting the last inequality, we obtain
\[
r_{2}-r \geq 
\delta(i_{0},j_{0},x,y) +
\delta(i_{1},j_{1},x_{1},y_{1}) -
\delta(i_{0},j_{0},x_{1},y_{1})-
\delta(i_{1},j_{1},x,y).
\]
Thereafter shifting terms we obtain
\[
\delta(i_{1},j_{1},x,y)-
\delta(i_{0},j_{0},x,y) 
\geq
\delta(i_{1},j_{1},x_{1},y_{1}) -
\delta(i_{0},j_{0},x_{1},y_{1}) + r - r_{2},
\]
and using the fact that 
\(s\geq r_{2}\)
it follows that
\[
\delta(i_{1},j_{1},x,y)-
\delta(i_{0},j_{0},x,y) 
\geq
\delta(i_{1},j_{1},x_{1},y_{1}) -
\delta(i_{0},j_{0},x_{1},y_{1}) + r - s = r -\Delta -s 
\]
from which it follows that 
\[(x,y)\not \in H[(i_{0},j_{0}),(i_{1},j_{1}),r-\Delta -s]\]
which proves (a).

Finally (b) follows from (a) by taking \(s=0\). 
Q.E.D.
\newline
{\bf Remark.} Part (b) of Theorem 4 can perhaps be considered 
as the main result
of the paper. From this result we  note that the larger the value 
\(r=\delta(i_{0},j_{0},x_{1},y_{1})-
\alpha(i_{0},j_{0})-
\beta(x_{1},y_{1})\) is, and the smaller the difference
\(\Delta = 
\delta(i_{0},j_{0},x_{1},y_{1}) -
\delta(i_{1},j_{1},x_{1},y_{1})\) is, the larger part of the plane
will be 
part of the exclusion set. The reason, that
we have included part (c) of the theorem, is that there are occasions,
when it can be simpler to check the size of the hyperbolic set 
\(H[(i_{0},j_{0}),(x_{1},y_{1}),
r-\delta(i_{0},j_{0},x_{1},y_{1})-s]\)
than it
is to check the size of the hyperbolic set
\(
H[(i_{0},j_{0}),(i_{1},j_{1}),
r-\Delta-s].\)
\section{An auxiliary lemma.} 
In this section we shall prove an elementary 
lemma which relates cones and hyperbolas in the plane. We shall rely on
this lemma when formulating and proving the results in the next section.

Let
\(\delta(\cdot,\cdot,\cdot,\cdot)\) denote the Euclidean metric.

\begin{lem} Let \(a\) and \(b\) be real numbers such that
\(a>0\) and \(b<2a\), and consider the following set:
\[ HYP = \{(x,y)\in R^{2}:  \delta(x,y,-a,0)-
\delta(x,y,a,0) > b\}.\]
Define \(H_{0}\) and \(H_{0}^{\prime}\) by 
\[ H_{0}=\{(x,y)\in R^{2}: x \geq 0\},\]
\[ H_{0}^{\prime}=\{(x,y)\in R^{2}: x > 0\},\]
and for \(c>0\) define \(CONE[x_{0},y_{0},c]\) by
\[CONE[x_{0},y_{0},c]=\{(x,y)\in R^{2}: x > x_{0}, 
|(y-y_{0})/(x-x_{0})|\leq c\}\cup (x_{0},y_{0}).\]
Then \newline
(a): if \(b <0\) then \(H_{0}\subset HYP\),
\newline
(b): if \(b =0\) then \(H_{0}^{\prime}= HYP\),
\newline
(c): if \(0<b<2a\) and \(c\leq \sqrt{(4a^{2}-b^{2})/b^{2}}\)
then  \(CONE[x,y,c]\subset HYP\) if \((x,y)\in HYP\).
\end{lem}
{\bf Proof.} Suppose that \((x,y)\in H_{0}\). Then  
\(\delta(x,y,-a,0)-
\delta(x,y,a,0)\geq 0\),
 and hence if
\(b < 0\), \((x,y)\in HYP\). This proves (a). 

Suppose that \((x,y)\in H_{0}^{\prime}\). Then  
\(\delta(x,y,-a,0) -
\delta(x,y,a,0) > 0\) and conversely if 
\(\delta(x,y,-a,0) -
\delta(x,y,a,0) > 0\) then \(x>0\).
Hence if
\(b = 0\), \(H_{0}^{\prime}=HYP\). This proves (b). 

The proof of
(c) is more complicated. Let \(b>0\). Consider the equation
\[\delta(x,y,-a,0)-
\delta(x,y,a,0)= b.
\]
Taking the square of each side we obtain the equality
\[(x-a)^{2} + y^{2} + (x+a)^{2} + y^{2} - 
 2\sqrt{((x-a)^{2} + y^{2})((x+a)^{2} + y^{2})} = b^{2}.\]
Moving the term \(b^{2}\) to the left hand side and 
moving the square root to the right hand side and making some
simplifications we obtain the equation:
\[2x^{2} + 2y^{2} + 2a^{2} - b^{2} = 
 2\sqrt{((x-a)^{2} + y^{2})((x+a)^{2} + y^{2})}.\]

Then again taking 
the square of both sides we obtain 
\[4x^{4} + 4y^{4} + (2a^{2} - b^{2})^{2} 
+8x^{2}y^{2} + 8x^{2}a^{2} -4x^{2}b^{2}+
8y^{2}a^{2} -4y^{2}b^{2}=\]
\[4((x-a)^{2} + y^{2})((x+a)^{2} + y^{2}),\]
and the right hand side of this equality can be simplified to
\[
4x^{4} + 4a^{4} - 8x^{2}a^{2} +4y^{4} + 
8x^{2}y^{2} + 8y^{2}a^{2}.
\] 
Eliminating common terms we obtain the equation
\[
4a^{4} + b^{4}-4a^{2}b^{2} +8x^{2}a^{2}
-4x^{2}b^{2} - 4y^{2}b^{2}
=
4a^{4} - 8x^{2}a^{2},
\]
and then, making further elimination and moving the terms 
containing \(x\)-factors and \(y\)-factors to the left hand side and the other
terms to the right hand side, we obtain the equation:
\[
16x^{2}a^{2}- 4x^{2}b^{2}-4y^{2}b^{2}=-b^{4}+4a^{2}b^{2}.\]
Finally, changing signs and dividing
all terms by 
\(4b^{2}\),
this equation can be written in a more familiar form as

\[ y^{2} = (x^{2}-b^{2}/4)((4a^{2}/b^{2})-1).
\]

This last equation determines a double sided hyperbola, and from the
equation we also conclude that the lines 
\[ y= \sqrt{(4a^{2}-b^{2})/b^{2}}x\]
and
\[ y= -\sqrt{(4a^{2}-b^{2})/b^{2}}x\]
are the asymptotic lines of the hyperbola. Since 
\(c\leq \sqrt{(4a^{2}-b^{2})/b^{2}}\),
it is clear from well-known properties of the hyperbola, that
the cone \(CONE[x,y,c]\) is a subset of \(HYP\) as soon as 
\((x,y)\in HYP\). QED.

\section{Stopping criteria for the Euclidean distance.} 
We shall now 
consider the problem of how to find new admissible arcs when the
underlying distance-function
\(\delta(i,j,x,y) \) 
is the Euclidean distance, i.e.  
\begin{equation}
\delta(i,j,x,y) = \sqrt{(i-x)^{2}+(j-y)^{2}}.
\end{equation} 
We shall
first introduce
some further terminology and notations.
We define
\[
    NE[i,j] =\{(x,y)\in R^{2}: x\geq i, y \geq j\}
\]
and
\[LE[i,j]=
\{(x,y)\in R^{2}: y=j, x \geq i\}.\]
(We have used the letters 
\(NE\) and \(LE\) as abbreviations of 
northeast and "line east").

Now let 
\((i_{0},j_{0})\) 
be a pixel in \(K_{1}\).
What we are interested in, is to find all pixels in
\(K_{2}\), which are admissible with respect to 
\((i_{0},j_{0})\). The purpose of this section is to present and prove some
results by which one can reduce the search time for finding the admissible
arcs with respect to
\((i_{0},j_{0})\), 
which also are {\bf \(NE\)} of 
\((i_{0},j_{0})\). By symmetry it is then easy to reformulate the results so
that they can be applied, when looking for pixels 
located \(NW\), \(SE\) or \(SW\)
of
\((i_{0},j_{0})\) 
which are admissible with respect to
\((i_{0},j_{0})\).

In case we have digital images defined on a
grid structure, the general algorithmic procedure is essentially as 
follows. First 
check all pixels on the 
line
\(LE[i_{0},j_{0}]\),
then on the line
\(LE[i_{0},j_{0}+1]\),
then on the line
\(LE[i_{0},j_{0}+2]\),
etcetera.

In order to reduce the search time there are (at least) two ways 
one can accomplish this. Firstly, for each line 
\(LE[i_{0},j_{0}+k], k=0,1,2,..\), 
 one can find a stopping
criterion which implies that one need not to check pixels further away
from 
\((i_{0},j_{0})\) on that line. Secondly, one can find 
a stopping criterion which implies that one does not have to check 
{\em any}
of the lines \(LE[i_{0},j_{0}+k]\)
above a certain value of \(k\).

In this section we shall present two theorems. The first gives  
stopping criteria along the lines
\(LE[i_{0},j_{0}+k], k>0\). 
The second gives stopping criteria
along the line \(\{(x,y):x=i_{0}, y\geq j_{0}\}\).
We formulate the theorems in such a way that they also
can be used when looking
for a quantity of the form  
\[ \min\{\delta(i,j,x,y)-\alpha(i,j)-\beta(x,y):
(i,j)\in A,\;\;(x,y)\in B\}\]
where \(A\subset K_{1}\) and  \(B\subset K_{2}\).
 
However before we state and prove our theorems we shall state
and prove a simple stopping criterion for pixels on the
line \(LE[i_{0},j_{0}]\). 
 
\begin{path} Suppose that the underlying distance-function is defined by (22).
Let 
\(\{\alpha(i,j), \beta(x,y)\}\) 
be a feasible set of dual variables such that 
each pixel \((i,j)\) in \(K_{1}\), and 
each pixel \((x,y)\) in \(K_{2}\), have an admissible arc.
Let \((i_{0},j_{0})\) belong to \(K_{1}\) and
\((x_{1},y_{1})\) belong to \(K_{2}\).
Suppose further that \((x_{1},y_{1})\) is such that
\(x_{1} \geq i_{0}\) and \(y_{1}=j_{0}\), and that
\((x_{1},y_{1})\) is not an admissible pixel with respect to 
\((i_{0},j_{0})\). Then 
\(LE[x_{1},y_{1}]\cap K_{2} \subset E[(i_{0},j_{0}),K_{2}]\). (Recall that
\(E[(i,j),K_{2}]\) denotes the exclusion set with respect to \((i,j)\)).
\end{path}
{\bf Proof.} Our proof follows the same line as our previous proofs.
Assume that \((x,y)\in LE[x_{1},y_{1}]\cap K_{2}\) and assume also that
\((x,y)\) is admissible with respect to \((i_{0},j_{0})\). We then have 
the following two relations:

\[ \delta(i_{0},j_{0},x,y)=\alpha(i_{0},j_{0})+\beta(x,y)\]
and
\[\delta(i_{0},j_{0},x_{1},y_{1}) > \alpha(i_{0},j_{0}) + 
\beta(x_{1},y_{1}).\]
By subtracting the inequality from the equality we obtain
\[\delta(i_{0},j_{0},x,y)-
\delta(i_{0},j_{0},x_{1},y_{1}) < \beta(x,y) -\beta(x_{1},y_{1}).\]
Since \(y=y_{1}=j_{0}\) and \(x\geq x_{1}\) it follows that
\[\delta(i_{0},j_{0},x,y)-
\delta(i_{0},j_{0},x_{1},y_{1})=x-x_{1}\] 
and consequently it follows that
\begin{equation}
x-x_{1} < \beta(x,y) -\beta(x_{1},y).
\end{equation}
But since \(\delta(\cdot,\cdot,\cdot,\cdot)\) is a metric
and since we have assumed that all pixels have admissible arcs,
 it follows
from Proposition 2 that
\(\beta(x,y) -\beta(x_{1},y) \leq \delta(x,y,x_{1},y)\)
and since  \(\delta(\cdot,\cdot,\cdot,\cdot)\) is the Euclidean metric
it follows that 
\[\delta(x,y,x_{1},y)=x-x_{1}.\] 
Hence
\(\beta(x,y) -\beta(x_{1},y) \leq x-x_{1}\)
which combined with (23) implies that
\(x-x_{1} < x-x_{1}\), by which we have reached 
a contradiction. QED.

This simple proposition gives rise to a 
stopping criterion when searching for admissible arcs on the line
\(y=j_{0}\).
Our next aim is to present a result which will be 
useful as a
stopping criterion when searching for admissible pixels 
in the sets
\(LE[i_{0},j_{0}+k], \;k=1,2,...\) .

\begin{thm}
Suppose that the underlying distance-function is defined by
\begin{equation}
\delta(i,j,x,y)=\sqrt{(i-x)^{2}+(j-y)^{2}}.
\end{equation}
Let 
\(\{\alpha(i,j), \beta(x,y)\}\) 
be a feasible set of dual variables such that 
each pixel \((i,j)\) in \(K_{1}\), and 
each pixel \((x,y)\) in \(K_{2}\), belong to an admissible arc.
Let 
\((i_{0},j_{0})\) 
belong to \(K_{1}\). Let
\((x_{1},y_{1})\) be a pixel in \(K_{2}\), 
NE of \((i_{0},j_{0})\), and different
from \((i_{0},j_{0})\), 
which is not an admissible pixel with respect to 
\((i_{0},j_{0})\).
Let instead 
\((i_{1},j_{1})\)
be an admissible pixel with respect to 
\((x_{1},y_{1})\), and
assume also that
\((i_{1},j_{1})\) is \(NE\) of 
\((i_{0},j_{0})\) and that 
\[i_{1}>i_{0}.\]
Let \(r\) be defined by
\begin{equation}
r=\delta(i_{0},j_{0},x_{1},y_{1})-
\alpha(i_{0},j_{0})-\beta(x_{1},y_{1}),
\end{equation}
let \(a\) be defined by
\begin{equation}
a=\delta(i_{0},j_{0},i_{1},j_{1})/2,
\end{equation}
let \(s\) be a real number such that 
\[ 
0\leq s < r,
\]
let \(b\) be defined by
\begin{equation}
b=\delta(i_{0},j_{0},x_{1},y_{1})-\delta(i_{1},j_{1},x_{1},y_{1}) - r + s,
\end{equation}
and let us also assume that \(b>0\). (The case when \(b\leq 0\) will 
be considered in the next theorem.)

Then,
if 
\[\sqrt{(4a^{2}-b^{2})/b^{2}}\geq (j_{1}-j_{0})/(i_{1}-i_{0}),\]
then 
\[LE[x_{1},y_{1}] \cap K_{2} \subset L[(i_{0},j_{0}),K_{2},s].\]
\end{thm}
{\bf Proof.} Let us consider the hyperbolic set
\(H[(i_{0},j_{0}), (i_{1},j_{1}),r-\Delta-s]\)
where \(r\) is defined by (25) and where \(\Delta\) is defined by  
\(\Delta = 
\delta(i_{0},j_{0},x_{1},y_{1}) -
\delta(i_{1},j_{1},x_{1},y_{1})\).
Let the set \(H\) be defined by
\begin{equation}
H=\{(x,y)\in R^{2}: \delta(x,y,i_{0},j_{0})-
\delta(x,y,i_{1},j_{1})>b\}
\end{equation}
where \(\delta(\cdot,\cdot,\cdot,\cdot)\) is defined by (24) and
the number \(b\) satisfies (27). From the definition of 
\(H[(i_{0},j_{0}),(i_{1},j_{1}),r-\Delta-s]\)
it is clear that we have
\[ H = H[(i_{0},j_{0}),(i_{1},j_{1}),r-\Delta-s].
\]
Let us also observe that
since \(0\leq s < r\) it follows that
\[\delta(x_{1},y_{1},i_{0},j_{0})-
\delta(x_{1},y_{1},i_{1},j_{1}) > \]
\[\delta(x_{1},y_{1},i_{0},j_{0})-
\delta(x_{1},y_{1},i_{1},j_{1})-r+s = b, \]
and therefore 
\((x_{1},y_{1})\) belongs to \(H\).

Since we have assumed that \(b>0\) and that \(i_{1}> i_{0}\),
the boundary 
\(\delta H\)
of \(H\) will 
be a hyperbola, and the vertex of 
\(\delta H\) will be along the line between the two points
\((i_{0},j_{0})\)
and
\((i_{1},j_{1})\),
and closer to 
\((i_{1},j_{1})\)
than to
\((i_{0},j_{0})\).
The axis will be in the \(NE\)-direction from
the vertex, since 
\((i_{1},j_{1})\)
is assumed to be \(NE\) of
\((i_{0},j_{0})\). 

In order to guarantee that the set \(LE[x_{1},y_{1}]\) 
belongs to \(H\) it is sufficient that the cone
\(C_{0}\) 
defined by 
\[C_{0}= \{(x,y) \in R^{2}: x> x_{1}, |(y-y_{1})/(x-x_{1})| \leq 
(j_{1}-j_{0})/(i_{1}-i_{0})\} \]
is contained in 
\(H\).
This happens if the angle
between the axis of the hyperbola
\(\delta H\) 
and the asymptotes of
\(\delta H\) is such that the tangent of the angle
is larger than or equal to 
\((j_{1}-j_{0})/(i_{1}-i_{0})\).
We call this angle 
\(\theta\).

Now by using the estimates obtained in Lemma 1 we conclude that if
the numbers \(a\) and \(b\) as defined by (26) and (27) satisfy
\begin{equation}
\sqrt{(4a^{2}-b^{2})/b^{2}}\geq (j_{1}-j_{0})/(i_{1}-i_{0}),
\end{equation}
then the angle 
\(\theta\) 
between the axis and the asymptotics is so large that
any cone 
\[
\{(x,y) \in R^{2}: x> x^{\prime}, 
|(y-y^{\prime})/(x-x^{\prime})| \leq 
(j_{1}-j_{0})/(i_{1}-i_{0})\} \]
is a subset of  \(H\)
if 
\((x^{\prime},y^{\prime})\) belongs to \(H\). Since we have already
proved above, that 
\((x_{1},y_{1})\) always 
belongs to \(H\), 
it follows that the cone \(C_{0}\) is a subset 
of \(H\).
Since evidently 
\(LE[x_{1},y_{1}]\) is a subset of \(C_{0}\)
it now follows from part (a) of 
Theorem 4 
that
\[LE[x_{1},y_{1}]\cap K_{2}\subset 
L[(i_{0},j_{0}),K_{2},s].\]
QED.

Proposition 4 and Theorem 5 are useful for giving stopping criteria
along lines parallel to the line \(y=j_{0}\). Our next theorem gives rise
to stopping criteria along the line \(x=i_{0}\). The basic assumptions
of the next theorem are similar to those of Theorem 5.

\begin{thm}
Again, suppose that the underlying distance-function is defined by
\[
\delta(i,j,x,y)=\sqrt{(i-x)^{2}+(j-y)^{2}},
\]  
and let as usual
\(\{\alpha(i,j), \beta(x,y)\}\) 
be a feasible set of dual variables such that 
each pixel \((i,j)\) in \(K_{1}\), and 
each pixel \((x,y)\) in \(K_{2}\), belong to an admissible arc.
Let 
\((i_{0},j_{0})\) 
belong to \(K_{1}\). Let
\((x_{1},y_{1})\) be a pixel in \(K_{2}\), 
\(NE\) of \((i_{0},j_{0})\), and different
from \((i_{0},j_{0})\), 
which is not an admissible pixel related to 
\((i_{0},j_{0})\).
Let instead 
\((i_{1},j_{1})\)
be an admissible pixel related to 
\((x_{1},y_{1})\) and assume that also 
\((i_{1},j_{1})\) is \(NE\) of 
\((i_{0},j_{0})\) (\(i_{1}\) not necessarily larger than \(i_{0}\)).
Again, let \(r\) be defined by
\[
r=\delta(i_{0},j_{0},x_{1},y_{1})-
\alpha(i_{0},j_{0})-\beta(x_{1},y_{1}),
\]
let \(a\) be defined by
\[
a=\delta(i_{0},j_{0},i_{1},j_{1})/2,
\]
let \(s\) be a real number such that 
\(\;\;
0\leq s < r,\;\;
\)
and let
\(b\) be defined by
\[
b=\delta(i_{0},j_{0},x_{1},y_{1})-\delta(i_{1},j_{1},x_{1},y_{1}) - r + s.
\]
(a): Suppose \(b \leq 0\) and that \(j_{1}>j_{0}\). Define
\begin{equation}
y_{2}=(j_{0}+j_{1})/2+(i_{1}-i_{0})^{2}/2(j_{1}-j_{0}).
\end{equation}
Then, if \(y >y_{2}\),
\[
    NE[i_{0},y] \cap K_{2} \subset L[(i_{0},j_{0}),K_{2},s].
\]

Next suppose that \(b>0\), that \(j_{1}>j_{0}\)
and also that \(i_{1}>i_{0}\). 
Suppose also  
\begin{equation}
j_{1}-j_{0}\leq i_{1}-i_{0}.
\end{equation}
Then  \newline
(b): if 
\[
\sqrt{(4a^{2}-b^{2})/b^{2}}\geq (i_{1}-i_{0})/(j_{1}-j_{0})
\]
\newline
then
\[
    NE[x_{1},y_{1}] \cap K_{2} \subset L[(i_{0},j_{0}),K_{2},s].
\]

Next, suppose instead of (31) 
 that 
\(
j_{1}-j_{0} > i_{1}-i_{0}.
\)
\newline
Then 
\newline
(c): if
\[\sqrt{(4a^{2}-b^{2})/b^{2}}\geq (j_{1}-j_{0})/(i_{1}-i_{0})\]
then 
  \[
    NE[x_{1},y_{1}] \cap K_{2} \subset L[(i_{0},j_{0}),K_{2},s],
\]
and 
\[
    NE[i_{0},y] \cap K_{2} \subset L[(i_{0},j_{0}),K_{2},s] 
\]
if \(y\geq y_{3},\)
where the point \(y_{3}\) is defined by 
\begin{equation}
y_{3}=
j_{1}+2(i_{1}-i_{0})^{2}(j_{1}-j_{0})/
((j_{1}-j_{0})^{2}-(i_{1}-i_{0})^{2}).
\end{equation}
\end{thm}
{\bf Proof.} Let us again consider the hyperbolic set
\(H[(i_{0},j_{0}), (i_{1},j_{1}),r-\Delta-s]\)
where \(r\) is defined by (25), and where \(\Delta\) is defined by  
\(\Delta = 
\delta(i_{0},j_{0},x_{1},y_{1}) -
\delta(i_{1},j_{1},x_{1},y_{1})\).
For each \(b\) define the set \(H(b)\) by
\[
H(b)=\{(x,y)\in R^{2}: \delta(x,y,i_{0},j_{0})-
\delta(x,y,i_{1},j_{1})>b\}
\]
where \(\delta(\cdot,\cdot,\cdot,\cdot)\) is defined by (24) and
the number \(b\) satisfies (27). From the definition of 
\(H[(i_{0},j_{0}),(i_{1},j_{1}),r-\Delta-s]\)
it is clear that we have
\[ H(b) = H[(i_{0},j_{0}),(i_{1},j_{1}),r-\Delta-s].
\]
From the definition of \(H(b)\) it is also clear that
\[b_{1}<b_{2} \Rightarrow H(b_{2})\subset H(b_{1}).\]
As we showed in the proof of the previous theorem
we also know that
\((x_{1},y_{1})\) belongs to \(H(b)\).

To prove (a) let us first consider the case when \(b=0\). The boundary
of  \(H(b)\) will in this case be a straight line, \(L\) say, defined
by the equation
\[y-(j_{0}+j_{1})/2 = - 
((i_{1}-i_{0})/(j_{1}-j_{0}))(x-(i_{0}+i_{1})/2).
\]
If we consider this equation as a function of \(x\) it is clear,
since \(i_{1}\geq i_{0}\) and \(j_{1}>j_{0}\), that the function
is non-increasing, and therefore, since \((x_{1},y_{1})\in H(b)\),
it is clear that 
\begin{equation}
NE[x_{1},y_{1}]\subset H(b)
\end{equation}
and also
that  
\begin{equation}
NE[i_{0},y]\subset H(b)
\end{equation}
if \(y>y_{2}\) where \(y_{2}\) is such that \((i_{0},y_{2})\) are the
coordinates of the point where the line \(L\) cuts the line whose 
equation is \(x=i_{0}\). 
To determine the value of \(y_{2 }\) we just have
to insert the value  
\(x=i_{0}\) into the equation defining \(L\). We then find that 

\[ y = (j_{0}+j_{1})/2 + ((i_{1}-i_{0})^{2}/2(j_{1}-j_{0}),\]
and hence \(y_{2}\) satisfies (29). 
That the assertions of (a) hold, now follows from (33) and (34), and hence 
(a) is proved in case \(b=0\). But that the assertions of (a) are true
also in case \(b<0\) follows immediately from the fact that 
\(b_{1}<b_{2} \Rightarrow H(b_{2})\subset H(b_{1}).\) Hence part (a)
of the theorem is proved.

From now on we
denote \(H(b)\) by \(H\). 
Next assume that \(b>0\) and that \(i_{1}> i_{0}\). 
Then the boundary 
of \(H\), which we denote by \(\delta H\), will 
again be a hyperbola, the vertex of 
\(\delta H\) will be along the line between two points
\((i_{0},j_{0})\)
and
\((i_{1},j_{1})\)
and closer to 
\((i_{1},j_{1})\)
than to
\((i_{0},j_{0})\), and the axis will be in the \(NE\)-direction from
the vertex, since 
\((i_{1},j_{1})\)
is assumed to be \(NE\) of
\((i_{0},j_{0})\). (Recall that \(H\) is defined by (28)).

We now want to find conditions on the numbers \(a\) 
and \(b\) as defined by (26) and (27) such that
the hyperbola \(\delta H\) has an eccentricity which is so large
that the set \(NE[x_{1},y_{1}]\) is contained in \(H\). Since 
we have assumed that  
\(j_{1}- j_{0} \leq i_{1}- i_{0}\),
the crucial angle is the angle between 
the \(y\)-axis and the line between the points
\((i_{0},j_{0})\)
and
\((i_{1},j_{1})\) (that is the line along the axis of the
hyperbola). Since we have assumed 
that \(j_{1}>j_{0}\), this angle is less than \(\pi/2\) and 
the tangent of the angle
is clearly equal to  
\((i_{1}- i_{0})/ (j_{1}- j_{0})\).
By Lemma 1 we conclude that the hyperbola \(\delta H\)
will have a sufficiently large eccentricity if the numbers 
\(a\) 
and \(b\) will satisfy the equality
\[\sqrt{(4a^{2}-b^{2})/b^{2}}\geq (i_{1}-i_{0})/(j_{1}-j_{0}).\]
instead of (29).
From Theorem 4 it thus follows that 
\[
NE[x_{1},y_{1}] \cap K_{2} \subset L[(i_{0},j_{0}),K_{2},s]\] 
and thus 
(b) is proved.

It remains to prove (c). Thus assume that \(b>0\), that 
\(i_{1}>i_{0}\), and that
\(j_{1}- j_{0} > i_{1}- i_{0}\).
Again we consider the set \(H\) defined by (28), and again we want the 
eccentricity of the hyperbola  \(\delta H\) to be so large that
the set \(NE[x_{1},y_{1}]\) is a subset of \(H\). This time the 
crucial angle is the angle between the \(x\)-axis and the 
line through the
two points
\((i_{0},j_{0})\)
and
\((i_{1},j_{1})\) (that is the line along the axis of the
hyperbola). Since we have assumed 
that \(i_{1}>i_{0}\), this angle is less than \(\pi/2\), and 
the tangent of the angle
is clearly equal to  
\((j_{1}- j_{0})/ (i_{1}- i_{0})\).
Again, by applying Lemma 1, we conclude that the hyperbola \(\delta H\)
will have a sufficiently large eccentricity if the numbers 
\(a\) 
and \(b\) will satisfy the inequality (29).
Hence by Theorem 4 it follows that if (29) holds then
\[
NE[x_{1},y_{1}] \cap K_{2} \subset L[(i_{0},j_{0}),K_{2},s]\] 
and thus the first part of 
(c) is proved.

Moreover,  if 
\(a\) 
and \(b\) satisfy the inequality (29), then the line, which
starts at the point
\((i_{1},j_{1})\) and which makes an angle with the axis of the hyperbola
\(\delta H\) which is the same as 
the angle between the axis of the hyperbola and
the \(x\)-axis, will be contained in \(H\). This line will cut the
line \(x=i_{0}\) at some point \(y_{3}\).
From geometric considerations it is then clear
that 
\(NE[i_{0},y,]\subset H\)
if 
\(y\geq y_{3}\)
and hence by Theorem 4 we also have
\[NE[i_{0},y] \cap K_{2} \subset L[i_{0},j_{0},K_{2},s]\] 
if 
\(y\geq y_{3}\).

It remains to show that 
the number \(y_{3}\) is determined by the expression (32). 
The equation for the line we are considering 
can be written 
\[
y-j_{1}=k(x-i_{1}),\]
where the number \(k\) still has to be determined. If we denote 
by
\(\theta\)
the angle between the line through the two points
\((i_{0},j_{0})\)
and
\((i_{1},j_{1})\) (that is the line along the axis of the
hyperbola \(\delta H\)) and the line \(y=i_{0}\) , then 
\[
k=tan(2\theta)=2tan(\theta)/(1-tan^{2}(\theta)).\]
Now since \(tan(\theta)=
(j_{1}-j_{0})/(i_{1}-i_{0})\),
by inserting
\(x=i_{0}\) into the equation for the line ,
we find that
\[
y=j_{1}+2(j_{1}-j_{0})(i_{0}-i_{1})/
(1-((j_{1}-j_{0})/(i_{1}-i_{0}))^{2},\]
which after simplification is equal to the expression in (32 ).

Thereby also the second part of (c) is proved and thereby the 
proof 
is completed. Q.E.D.
\newline
{\bf Remark.} 
By using the theorems of this section, it seems likely, that 
the search time for finding all pixels which are 
admissible with respect to a given
pixel in most situations will have a 
computational complexity of order \(O(1)\),
and therefore the computational complexity for finding all admissible arcs
will be of order \(O(N)\) roughly, where though the constant in \(O(N)\)
may be fairly large.

\section{Each pixel belongs to an admissible arc.}
One of the basic assumptions we have made in most of the results 
proven above is that the set of dual feasible variables 
\(\{\alpha(i,j),\beta(x,y)\}\) 
is such that each pixel 
\((i,j)\in K_{1}\)
and each pixel
\((x,y)\in K_{2}\) belong to an admissible arc. That is, to each
\((i,j)\in K_{1}\) there exists a pixel 
\((x,y)\in K_{2}\) such that
\begin{equation}
\alpha(i,j)+ \beta(x,y)=\delta(i,j,x,y), 
\end{equation}
and to each 
\((x,y)\in K_{2}\) there exists a pixel
\((i,j)\in K_{1}\) such that (35) holds. 
 We shall now
show that this is a valid assumption independently of the choice of
the underlying distance-function.

Let us first state the following proposition.
\begin{path} 
Define the set of dual variables 
\(\{\alpha(i,j),\beta(x,y)\}\) as follows:
\begin{equation}
 \alpha(i,j)=\min\{\delta(i,j,x,y):(x,y)\in K_{2}\}
\end{equation}
\begin{equation}
\beta(x,y)=\min\{\delta(i,j,x,y)-\alpha(i,j):(i,j)\in K_{1}\}.
\end{equation}
Then the set of dual variables is feasible, and, moreover, to 
each pixel 
\((i,j)\in K_{1}\)
and to each pixel
\((x,y)\in K_{2}\) 
there exists an admissible arc.
\end{path}
{\bf Proof.} The proposition is intuitively obvious.
A formal proof can read as follows.
>From equation (37) it is clear that to each 
\((x,y)\in K_{2}\) 
there exists a pixel
\((i,j)\in K_{1}\) 
such that (35) holds. Moreover, from (37) it also follows
that the set of dual variables is a feasible set. 
Now suppose that there
exists a pixel 
\((i,j)\in K_{1}\) 
such that 
\[
\alpha(i,j)+ \beta(x,y)<\delta(i,j,x,y)
\]
for all
\((x,y)\in K_{2}\). Since 
\(\alpha(i,j)=\delta(i,j,x,y)\)
for some
\((x,y)\in K_{2}\), it follows that for this choice of pixel,
\(\beta(x,y)\) must be negative. But this is
not possible because of (37) and the fact that 
\(\delta(i,j,x,y)\geq \alpha(i,j)\)
because of (36). Q.E.D.

We shall now prove:
\begin{lem}Let the set 
\(\{\alpha(i,j),\beta(x,y)\}\) be a set of feasible dual variables
obtained by the primal-dual algorithm when initially the set of 
dual variables are defined as in Proposition 5. Then 
each pixel 
\((i,j)\in K_{1}\) 
and each pixel 
\((x,y)\in K_{2}\) will 
have an admissible arc.
\end{lem} 
{\bf Remark.} In case one is familiar with
the primal-dual algorithm, the truth of the lemma is intuitively clear. 
We include
a formal proof for the sake of completeness.
\newline
{\bf Proof.} To give a formal proof of the lemma we need to describe
the primal-dual algorithm partly again. As pointed out in
Section 4, one step in the primal-dual
algorithm is called the labeling procedure, and below we shall 
describe it in more detail. Let us first introduce two
further notions. Thus let 
\(h(i,j,x,y)\) 
be a transition plan from the image 
\(P=\{p(i,j):\;(i,j)\in K_{1}\}\) 
to the image
\(Q=\{q(x,y):\;(x,y)\in K_{2}\}\). (See Section 2 for the
definition of a transportation plan.)
We call a pixel \((i,j)\in K_{1}\)  \emph{deficient} if
the transportation plan
\(h(i,j,x,y)\) is such that
\begin{equation}
\sum_{(x,y) \in K_{2}} h(i,j,x,y) < p(i,j). 
\end{equation}
If instead \((i,j)\in K_{1}\) is such that the left hand side
of (38) is \emph{equal} to the right hand side of (38) then we
say that \((i,j)\) is \emph{full}. 

We denote the subset of \(K_{1}\)
consisting of all deficient pixels by \(D_{1}\) and the subset 
consisting of all full pixels by \(F_{1}\).

The labeling procedure is as follows. 
One starts the labeling procedure by labeling
all deficient elements
in \(K_{1}\) (that is the set \(D_{1}\)). Then,
whenever a pixel 
\((i,j)\in K_{1}\) is labeled, then all pixels 
\((x,y)\in K_{2}\) which are admissible with respect to
\((i,j)\) and are not yet labeled, are labeled, and
whenever a pixel 
\((x,y)\in K_{2}\) is labeled, then all pixels \((i,j)\in K_{1}\),
 which are admissible with respect to \((x,y)\), which are not yet 
labeled and for which \(h(i,j,x,y)>0\), are labeled.

At the end of a labeling procedure, depending on the labeling, 
one can either go to the so called {\em flow change routine} or to the
so called 
{\em dual solution change routine}. If one goes to the flow change routine
the transportation plan is redefined after which one can redo the
labeling. Sooner or later one has reached a state when it is not possible to 
find a better transportation plan for the given set of admissible arcs
induced by the present set of dual variables. One therefore has
to redefine the dual variables.

When the labeling procedure ends, a number of pixels in both 
\(K_{1}\)
and
\(K_{2}\) are labeled. 
Let  
\(L_{1}\)
denote the set of labeled pixels in
\(K_{1}\),
let
\(U_{1}\)
denote the set of unlabeled pixels in
\(K_{1}\),
let  
\(L_{2}\)
denote the set of labeled pixels in
\(K_{2}\),
and let
\(U_{2}\)
denote the set of unlabeled pixels in
\(K_{2}\).

When the dual variables are changed, one uses the following quantity:
\begin{equation}
\Theta = \max \{\delta(x,y,i,j)-\alpha(i,j)-\beta(x,y):
(i,j)\in L_{1}, (x,y)\in U_{2}\}
\end{equation}
which one can prove will be \( > 0\) (unless the transportation plan
is complete). The way the dual variables are changed is as follows:
\[\alpha_{new}(i,j) = \alpha_{old}(i,j)+\Theta, \;\forall (i,j) \in L_{1}\]
\[\alpha_{new}(i,j)=\alpha_{old}(i,j), \;\forall (i,j) \in U_{1}\]
\[\beta_{new}(x,y)=\beta_{old}(x,y)-\Theta, \;\forall (x,y) \in L_{2}\]
\[\beta_{new}(x,y)=\beta_{old}(x,y), \;\forall (x,y) \in U_{2}.\]
From the way the new set of dual
variables are defined the following conclusion follows immediately. 
\begin{path}
After we have updated the dual variables the only way 
by which an arc 
\(\{(i,j),(x,y)\}\)
can
be unadmissible when it previously was admissible is
 if \((i,j)\in U_{1}\) and \((x,y)\in L_{2}\).
\end{path}
Now let us thus prove that the new set of dual variables
gives rise to a set of admissible arcs such that
to each pixel in 
\(K_{1}\) 
and to each pixel in
\(K_{2}\) there exists an admissible arc.

Our induction assumption is that the old set of dual variables has this 
property. Note that initially this is true because of Proposition 5.
Now let us first prove that 
to each pixel \((i,j)\in 
K_{1}\) 
there  exists a pixel in
\(K_{2}\) which is an admissible arc.
From the induction hypothesis and Proposition 6 above, the only case
we have to consider is the case when \((i,j)\) is unlabeled. But if
\((i,j)\) is unlabeled it must be full. Next let \(G(i,j)\) be
the set of pixels \((x,y)\in K_{2}\) for which \(\{(i,j),(x,y)\}\)
is an admissible arc under the old set of dual variables. If there
exists an element \((x,y)\in G(i,j)\) which is not labeled then 
we have found an arc which will also be admissible under the new
set of dual variables. Otherwise all pixels in \(G(i,j)\) must be
labeled. Since the pixel \((i,j)\) is full it follows that there
must exist a pixel \((x,y)\in G(i,j)\) for which \(h(i,j,x,y)>0\).
But then because of the way the labeling procedure works the pixel
\((i,j)\) would have been labeled. Hence if there exists 
a pixel \((i,j)\in K_{1}\) 
which is not labeled
there must exist a pixel \((x,y)\in G(i,j)\) which is not labeled 
and consequently this arc will be admissible also after the
dual variables are updated. 

We have now proved that to every pixel
\((i,j)\in K_{1}\) there exists a pixel 
\((x,y)\in K_{2}\) such that the arc \(\{(i,j),(x,y)\}\) is
admissible under the new set of dual variables. It remains to prove that
each pixel 
\((x,y)\in K_{2}\) 
has an admissible arc under the new set of dual variables.
But this is somewhat easier.
By Proposition 6, a necessary condition for an arc 
\(\{(i,j),(x,y)\}\)
to
be unadmissible when it previously was admissible is
that \((x,y)\) is labeled. But from the labeling procedure, it then
must exist an admissible pixel \((i,j)\) with respect to \((x,y)\)
which also is labeled, and then because of
Proposition 6 the arc \(\{(i,j),(x,y)\}\) will
be admissible also under the new set of dual variables. Thereby the induction
step is verified and the proof is completed. Q.E.D

\section{Reducing the computation time for determining the quantity
by which the dual variables are changed.} 
An important step in the primal-dual algorithm for the 
transportation problem is
to determine the quantity \(\Theta\) as defined by (39). In general
the computational complexity for this part of the algorithm is
\(O(N^{2})\), where though 
\(O(N^{2})\), so to speak, has
a small constant.
 
In case one has integer storages and demands and 
integer-valued cost-function then  
one can always take \(\Theta\) equal to \(1\) and this brings the computation 
time down to zero.
This strategy 
works also quite well in case we are working with digital images. 

In case 
we are dealing with the transportation problem in the plane and have 
{\em the Euclidean distance} 
as cost-function and the positions of the
sources and sinks are essentially random, then it is necessary to 
compute the quantity \(\Theta\) defined by (39). 
In this section we shall formulate a result by which one can 
reduce the computation time to determine the quantity \(\Theta\). 

We first prepare for our result.
Let as always \(K_{1}=\{(i(n),j(n): n=1,2,...,N\}\)
and \(K_{2}=\{(x(m),y(m): m=1,2,...,M\}\) be two sets in \(R^{2}\), let 
\(P=\{p(i,j):(i,j)\in K_{1}\}\) and
\(Q=\{q(x,y):(x,y)\in K_{2}\}\) be two images such that
\(\sum_{K_{1}}p(i,j) = \sum_{K_{2}}q(x,y)\).
Let the underlying distance-function be defined by
\begin{equation}
\delta(i,j,x,y)=\sqrt{(i-x)^{2}+(j-y)^{2}}.
\end{equation}
Let 
\(\{\alpha(i,j), \beta(x,y)\}\) 
be a feasible set of dual variables such that 
each pixel \((i,j)\) in \(K_{1}\), and 
each pixel \((x,y)\) in \(K_{2}\) belong to an admissible arc.
Let \(L \subset K_{1}\) and \(U \subset K_{2}\) be such that
\begin{equation}
\Theta = \min\{\delta(i,j,x,y)-\alpha(i,j)-\beta(x,y):
(i,j)\in L, (x,y)\in U\} \, > \,0.
\end{equation}
Our aim is now to determine \(\Theta\) and we do this in principal
by computing \[ \Theta_{(i,j)} = \min 
\{\delta(i,j,x,y)-\alpha(i,j)-\beta(x,y):(x,y)\in U\}\]
for each \((i,j)\in L\). 

Let the size of \(L\) 
be equal to \(N_{1}\), let the size of \(U\) be \(M_{1}\)
let 
\(\{(i(k),j(k)):k=1,2,...,N_{1}\}\) 
be a sequential list of the set \(L\), 
let
\(\{(x(k),y(k)):k=1,2,...,M_{1}\}\) be a sequential list of \(U\).
For \(n= 1,2,...,N_{1},\) and \(m = 1,2,...,M_{1},\)
let 
\[ \Theta(n) = \min\{ \Theta_{(i(k),j(k))}: 1\leq k\leq n\},\]
let
\[\Theta_{(i,j)}(m) = \min\{\delta(i,j,x(k),y(k))-\alpha(i,j)-
\beta(x(k),y(k)),1\leq k \leq m \}, \]
let 
\[\Theta(0,m)=\Theta_{(i(1),j(1))}(m),\]
set \[\Theta(0)=\Theta(0,1),\]
and define  
\begin{equation}
\Theta(n,m)= \min \{\Theta(n), \Theta_{(i(n+1),j(n+1))}(m)\}.
\end{equation}

\begin{thm} 
Let \(L =\{(i(k),j(k)):k=1,2,...,N_{1}\} \subset K_{1}\) and 
\newline
\(U=\{(x(k),y(k)):k=1,2,...,M_{1}\}
 \subset K_{2}\) be such that
(41) holds.
Let \(n\) and \(m\) be two integers such that
\(0\leq n \leq N_{1}-1\) and \(1\leq M_{1}-1\) and \(n+m \geq 1\).
Set \((i_{0},j_{0})=(i(n+1),j(n+1))\),
set \((x_{1},y_{1})=(x(m+1),y(m+1))\),
let \((i_{1},j_{1})\) be an admissible
pixel with respect to \((x_{1},y_{1})\) and assume that 
\((i_{1},j_{1})\) is NE of \((i_{0},j_{0})\). 
Let \(r\) be defined by
\[
r=\delta(i_{0},j_{0},x_{1},y_{1})-
\alpha(i_{0},j_{0})-\beta(x_{1},y_{1}),
\]
let \(a\) be defined by
\[
a=\delta(i_{0},j_{0},i_{1},j_{1})/2,
\]
and let \(b\) be defined by
\[
b=\delta(i_{0},j_{0},x_{1},y_{1})-\delta(i_{1},j_{1},x_{1},y_{1}) -r+
\Theta(n,m)
\]
where \(\Theta(n,m)\) is defined by (42).
\newline
Then:
\newline
(a): Suppose \(b \leq 0\) and that \(y_{1}\geq (j_{0}+j_{1})/2\). Then 
\(
    NE[x_{1},y_{1}] \cap K_{2} \subset L[(i_{0},j_{0}),K_{2},\Theta(n,m)].
\)
\newline
(b): Next suppose that \(b>0\) 
and suppose also that \(j_{1}>j_{0}\). Suppose also  
\begin{equation}
j_{1}-j_{0}\leq i_{1}-i_{0}
\end{equation}
and
\[\sqrt{(4a^{2}-b^{2})/b^{2}}\geq (i_{1}-i_{0})/(j_{1}-j_{0}).\]
Then
\[
    NE[x_{1},y_{1}] \cap K_{2} \subset L[(i_{0},j_{0}),K_{2},\Theta(n,m)].
\]
Suppose instead of (43)
 that 
\(
j_{1}-j_{0} > i_{1}-i_{0}>0.
\)
\newline
Then, if
\[\sqrt{(4a^{2}-b^{2})/b^{2}}\geq (j_{1}-j_{0})/(i_{1}-i_{0})\]
then 
  \[
    NE[x_{1},y_{1}] \cap K_{2} \subset L[(i_{0},j_{0}),K_{2},\Theta(n,m)].
\]
\end{thm}
{\bf Remark.}
Since the proof is very similar to the proof
of Theorem 6, we omit it.

\section{Organizing points along the northeast direction}
In order to use Theorems 6 and 7 efficiently, it is necessary to have
an initialization process such that 
one can find elements \(NE\) of a given element rapidly.
Thus, it is necessary that one can organize the elements in the sets 
\(K_{1}\)
and
\(K_{2}\) and also the set 
\(K_{1}\cup K_{2}\) in the northeast direction, - as well 
as the northwest, the southeast and the southwest directions. 
In case we are dealing with digital images on a rectangular grid then
the pixels are organized already from the start 
in such a way that this is no problem.
However, in the general case it is necessary to organize the elements of
the two given images in the \(NE-SW\) direction, and the \(NW-SE\) direction. 
In this section we shall describe
one algorithm by which this can be done. 
We shall only consider
the northeast direction since the other directions 
can be handled analogously.

Thus let \(K=\{(x(n),y(n)), n=1,2,...,N\}\) be a set of points in the
plane. What we want to do is to organize the points in
\(K\) so that we can apply Theorems 6 and 7 efficiently.
What we will create is a structure which is often called a {\em directed 
acyclic graph}.

We start by introducing an extra point \((x(0),
y(0))\) such that
\(x(0) < x(n),\; n=1,2,...,N\) 
and 
\(y(0) < y(n), \; n=1,2,...,N\).

Next let us define \(z(n)=x(n)+y(n), \; n=1,2,...,N\) and let us
first order the elements of \(K\) in such 
a way that
\[z(n) \leq z(n+1), n=1,2,...,N-1\]
and  \[z(n)=z(n+1)  \Rightarrow  x(n) < x(n+1).\] 
This way to order the elements of \(K\) implies of course that
\((x(n),y(n))\) can not be \(NE\) of
\((x(m),y(m))\) if  \(n<m\).

To each element \((x(n),y(n)), \;0\leq n \leq N  \) we shall associate
a number of related elements which will be called {\em parents}
and {\em children}. These notions are defined as follows:
Let 
\((x(n),y(n))\)
and
\((x(m),y(m))\) be two points in \(K\) with \(n<m\), and
suppose that 
\((x(n),y(n))\) is SW of 
\((x(m),y(m))\).
The point
\((x(n),y(n))\) will be a {\em parent} of 
\((x(m),y(m))\) if either 
1) it is the only element in \(K\) which is \(SW\) of 
\((x(m),y(m))\), or
2)
if there are other points in \(K\) which are  
\(SW\) of 
\((x(m),y(m))\) then
\((x(n),y(n))\) is {\em not SW} of any such point.
If 
\((x(n),y(n))\) is a parent of 
\((x(m),y(m))\), then we say that 
\((x(m),y(m))\) is a {\em child} of 
\((x(n),y(n))\).

In case we define parents and children as 
above, then, if we have found the parents of all 
points in a set \(K \in R^{2}\),
then we say that we have organized the set in the \(NE-SW\) direction. 

It is easy to construct examples of sets of size \(N\) for
which the number of parents will be \(N^{2}/4\). On the other hand it is
clear that one does not need more than \(N(N+1)/2\) checks to
find all parents. Therefore it is clear that
the computational complexity 
to organize a set in the 
\(NE-SW\) direction is of order \(O(N^{2})\). However in most situations
we believe that one can 
organize a set in the \(NE-SW\) direction substantially faster 
by applying the algorithm we shall now  describe.

Let us assume that the points \((x(i),y(i)), \; i=0,1,2,...,n\) have
been organized in the \(NE-SW\) direction. We shall now show
how to find the parents of the next point 
\((x(n+1),y(n+1))\).

Let \(m_{1}\) be defined as the largest index for which 
\[x(m_{1})= \max \{x(i):\, x(i) \leq x(n+1), \, y(i)\leq y(n+1),
\,0\leq i \leq n\}\]
and let
\(m_{2}\) be defined as the largest index for which 
\[y(m_{2})= \max \{y(j):\, x(j) \leq x(n+1), \, y(j)\leq y(n+1),
\; 0\leq i \leq n\}.\]

We now have two possibilities. Either 
\(m_{1}=m_{2}\) or 
\(m_{1}\neq m_{2}\). In the first case we are ready and 
\((x(m_{1}),y(m_{1}))\) is the only parent of 
\((x(n+1),y(n+1))\). 

If instead \(m_{1}\neq m_{2}\) then we have obtained two parents
to 
\newline
\((x(n+1),y(n+1))\) 
namely \((x(m_{1}),y(m_{1}))\) 
and
\((x(m_{2}),y(m_{2}))\), and it is possible that there
are further parents. In order to find these, we shall look 
for points \((x(k),y(k)) \in K\) satisfying
\[ x(m_{2})\leq x(k)< x(m_{1})\]
and
\[ y(m_{1})\leq y(k)< y(m_{2}).\]

We now define
\(m_{3}\) 
as the largest index for which 
\[x(m_{3})= \max \{x(i):\, x(m_{2})\leq x(i) < x(m_{1}), 
\, y(m_{1})\leq y(i) <y(m_{2}),
\,1\leq i \leq n\},\]
if any such index exists,
and we define 
\(m_{4}\) as the largest index for which 
\[y(m_{4})= \max \{y(j):\, x(m_{2})\leq x(j) < x(m_{1}), 
\,  y(m_{1})< y(j)< y(m_{2}),
\,1\leq j \leq n\}.\]
If \(m_{3}\) is not defined, then there is no further
parent. 
In case 
\(m_{3}=m_{4}\) 
then
there is exactly one more parent namely  
\((x(m_{3}),y(m_{3}))\).
 
Otherwise we have found two new parents 
namely  
\((x(m_{3}),y(m_{3}))\)
and 
\((x(m_{4}),y(m_{4}))\). In this case there may be further parents 
namely if there are elements \((x(k),y(k))\in K\) which satisfy
\[ x(m_{4})\leq x(k)< x(m_{3})\]
and
\[ y(m_{3})\leq y(k)< y(m_{4}).\]
To investigate whether there exist any points in this rectangle
we proceed in the same way as above. 

We have now briefly described an algorithm by which 
one can organize the elements 
of \(K_{1}, K_{2}\) and \(K_{1}\cup K_{2}\) in such a way that
it is easy to find the nearest elements in the 
\(NE-SW\) directions.
\newline
{\bf Remark.} 
As pointed out above it is easy to construct an example
of a set with \(N\) pixels which has \(N^{2}/4\) parents and children. 
Whether
this also requires that we need a storage for the structure of
order \(N^{2}\) is not 100 \%  sure, since in the more extreme
situations it is likely that many "parent sets" are very similar,
and therefore it is conceivable that these sets 
could be coded efficiently. In the
general random case we believe that the storage needed to
store the information about parents and children will 
be much less than \(N^{2}\) when \(N\) is large.
\[\]
\[\]
{\bf 14. SUMMARY.} 
\newline
In this paper we have proved some results,
which give rise to stopping criteria, when
searching for new admissible arcs, when 
using the primal-dual algorithm for solving the transportation 
problem in the plane, 
in case the underlying distance-function 
\(\delta(\cdot,\cdot,\cdot,\cdot)\) 
is defined either by
\(\delta(i,j,x,y)=(x-i)^{2}+(y-j)^{2}\)
or by
\(\delta(i,j,x,y)=\sqrt{(x-i)^{2}+(y-j)^{2}}.\)
We believe, that by 
using Theorems 5 and 6 of Section 10, and Theorem 7 of section 12, 
it will be possible to reduce the computational complexity
of the Euclidean transportation problem and the Euclidean assignment
problem to approximately \(O(N^{2})\) from approximately
\(O(N^{2.5})\), which is the best limit so far, as far as the author knows.
(See e.g. 
\cite{AV95}
and
\cite{Vai89}).

\[\]
{\bf ACKNOWLEDGEMENT.} It is a pleasure to thank Erik Ouchterlony
for stimulating discussions and valuable comments.

\end{document}